\documentclass[12pt]{amsart}
\usepackage{a4wide,amsfonts,amssymb,stmaryrd,amscd,amsmath,latexsym,amsbsy}
\numberwithin{equation}{section}
\theoremstyle{plain}
\newtheorem{thm}{Theorem}[section]

\newtheorem{prop}[thm]{Proposition}

\newtheorem{defi}[thm]{Definition}
\newtheorem{lem}[thm]{Lemma}
\newtheorem{cor}[thm]{Corollary}

\theoremstyle{remark}
\newtheorem{rema}[thm]{Remark}
\title{Askey-Wilson functions and quantum groups.} 
\author{Jasper V. Stokman}
\address{J.V. Stokman,
KdV Institute for Mathematics, Universiteit van Amsterdam,
Plantage Muidergracht 24, 1018 TV Amsterdam, The Netherlands.}
\email{jstokman@science.uva.nl}
\begin{document}

\begin{abstract}
Eigenfunctions of the Askey-Wilson 
second order $q$-difference operator for $0<q<1$ and
$|q|=1$ are constructed as formal 
matrix coefficients of the principal series 
representation of the quantized 
universal enveloping algebra $\mathcal{U}_q(\mathfrak{s}
\mathfrak{l}(2,\mathbb{C}))$. 
The eigenfunctions are in integral form and 
may be viewed as analogues of Euler's 
integral representation for Gauss' hypergeometric
series. We show that for $0<q<1$ the resulting
eigenfunction can be rewritten as a 
very-well-poised ${}_8\phi_7$-series, and reduces
for special parameter values to a natural elliptic analogue of 
the cosine kernel.
\end{abstract}

\maketitle
\begin{center}
{\it Dedicated to Mizan Rahman}
\end{center}

\vspace{.4cm}
{\bf Contents.}
\begin{enumerate}
\item[{\S 1}] Introduction.
\item[{\S 2}] Generalized gamma functions.
\item[{\S 3}] The Askey-Wilson second order 
difference operator and quantum groups.
\item[{\S 4}] The Askey-Wilson function for $0<q<1$.
\item[{\S 5}] The expansion formula and the elliptic cosine kernel.
\item[{\S 6}] The Askey-Wilson function for $|q|=1$.
\item[{\S 7}] Appendix.
\end{enumerate}
\section{Introduction}
The aim of this paper is to simplify 
the quantum group construction of explicit eigenfunctions of
the second order Askey-Wilson $q$-difference operator, 
and to extend the results to the interesting
and less well studied $|q|=1$ case.

The approach is based on the known fact from \cite{K}, \cite{NM} and \cite{Ko}
that the second order Askey-Wilson difference operator arises as radial part of
the quantum Casimir element of 
$\mathcal{U}_q(\mathfrak{s}\mathfrak{l}(2,\mathbb{C}))$ 
when the radial part is computed with respect
to Koornwinder's twisted primitive elements. Using this result,
we construct nonpolynomial eigenfunctions 
of the Askey-Wilson second order difference
operator as matrix coefficients of the principal 
series representation of 
$\mathcal{U}_q(\mathfrak{s}\mathfrak{l}(2,\mathbb{C}))$. 
The two cases $0<q<1$ and $|q|=1$ will be treated seperately. 
The theory for $0<q<1$ is related to the noncompact quantum group
$\mathcal{U}_q(\mathfrak{s}\mathfrak{u}(1,1))$, 
while for $|q|=1$ it is related to
the noncompact quantum group 
$\mathcal{U}_q(\mathfrak{s}\mathfrak{l}(2,\mathbb{R}))$.

This approach was considered
in \cite{KS} for $0<q<1$ using an explicit realization of the principal series 
representation on $l^2(\mathbb{Z})$. The resulting eigenfunction then appears
as a non-symmetric Poisson type kernel involving nonterminating
${}_2\phi_1$ series.
With the help of a highly nontrivial summation formula,
proved by Rahman in the appendix of \cite{KS} (see \cite{KR} for extensions), 
this eigenfunction was expressed as one of the
explicit ${}_8\phi_7$-solutions of the Askey-Wilson second order difference
operator from \cite{IR}.
This eigenfunction was called the Askey-Wilson function in \cite{KS2},
since it is a meromorphic continuation of the Askey-Wilson polynomial
in its degree. 
In this paper we start by reproving this result, now using an explicit
realization of the principal series representation of 
$\mathcal{U}_q(\mathfrak{s}\mathfrak{l}(2,\mathbb{C}))$ as difference operators
acting on analytic functions on the complex plane. Koornwinder's twisted
primitive element then acts as a first order difference operator,
hence eigenvectors are easily constructed (for the positive discrete
series, this was observed by Rosengren in \cite{R}). The corresponding matrix
coefficients lead to explicit integral 
representations for eigenfunctions of the 
Askey-Wilson second order difference operator. These matrix coefficients
can be rewritten as the explicit ${}_8\phi_7$-series representation of 
the Askey-Wilson function by a residue computation.

We also show that for a special choice of parameter
values, the Askey-Wilson function reduces to an elliptic analogue of the
cosine kernel. This is the analogue of the classical fact
that the Jacobi function reduces to the 
cosine kernel for special parameter values,
see e.g. \cite{Koexp}. We give two proofs, one proof uses an explicit expansion
formula of the Askey-Wilson function in Askey-Wilson polynomials from
\cite{St0}, the other proof uses Cherednik's Hecke algebra techniques
from \cite{C} and \cite{St}. 

In the second part of the paper we consider the quantum group techniques for 
$|q|=1$. In this case, the approach is similar to the construction
of quantum analogues of Whittaker vectors and Whittaker functions from
\cite{KLS}. The
role of $q$-shifted factorials, or equivalently $q$-gamma functions,
is now taken over by Ruijsenaars' \cite{Ru1} hyperbolic gamma function. 
The hyperbolic
gamma function is directly related to Barnes' double gamma function, 
as well as to Kurokawa's double sine function, see \cite{Ru2} and references
therein. The quantum group technique applied to this particular set-up
leads to an eigenfunction of the Askey-Wilson
second order difference operator for $|q|=1$, given explicitly
as an Euler type integral involving hyperbolic gamma functions.

The emphasis in this paper lies on exhibiting the similarities between
the $0<q<1$ case and the $|q|=1$ case as much as possible. 
Other approaches might very well lead to  
eigenfunctions for the Askey-Wilson second order difference operator 
for $|q|=1$ which are ``more optimal'', in the sense that they 
satisfy {\it two} Askey-Wilson type 
difference equations in the geometric parameter, one with respect to
base $q=\exp(2\pi i\tau)$, the other with respect to the ``modular
inverted'' base $\widetilde{q}=\exp(2\pi i/\tau)$, cf.
\cite{KLS} for $q$-Whittaker functions.
Such eigenfunctions are expected to be realized as 
matrix coefficients of the {\it modular double} of the 
quantum group
$\mathcal{U}_q(\mathfrak{s}\mathfrak{l}(2,\mathbb{C}))$
(a concept introduced by Faddeev in \cite{F}), 
and are expected to be closely related to 
Ruijsenaars' \cite{Ru2}, \cite{Ru3} R-function. The R-function is an 
eigenfunction of two Askey-Wilson type second order
difference operators in the geometric parameter, 
which is explicitly given
as a Barnes' type integral involving
hyperbolic gamma functions. 
I hope to return to these considerations in a future paper.

{\it Acknowledgments:}
It is a pleasure to
dedicate this paper to Mizan Rahman. His important contributions to the 
theory of basic hypergeometric series and, more concretely, his kind
help in the earlier stages of the research on the Askey-Wilson functions
in \cite{KS}, have played, and still play, an important
role in my research on Askey-Wilson functions. 

I am supported by the Royal Netherlands Academy
of Arts and Sciences (KNAW). 


\section{Generalized gamma functions}

In this section we discuss 
$q$-analogues of the gamma function for deformation parameter $q$ 
in the regions $0<|q|<1$ and $|q|=1$.

\subsection{The $q$-gamma function for $0<|q|<1$.}
Let $\tau$ be a fixed complex number in the upper half plane $\mathbb{H}$.
The corresponding deformation parameter 
$q=q_\tau=\exp(2\pi i\tau)$ has modulus 
less than one. We write $q^u=\exp(2\pi i\tau u)$ for $u\in\mathbb{C}$.

Let $b,b_j\in\mathbb{C}$ and $n\in\mathbb{Z}_{\geq 0}\cup\{\infty\}$.
The $q$-shifted factorial is defined by
\[\bigl(b;q\bigr)_n=\prod_{j=0}^{n-1}(1-bq^j),\qquad
\bigl(b_1,\ldots,b_m;q\bigr)_n=\prod_{j=1}^m\bigl(b_j;q\bigr)_n.
\]
In $q$-analysis the function
\[x\mapsto 
\frac{\bigl(q;q\bigr)_{\infty}}{\bigl(q^x;q\bigr)_{\infty}}(1-q)^{1-x}
\]
is known as the $q$-gamma function, 
see e.g. \cite{GR}.  For our purposes, it is more convenient 
to work with the function
\begin{equation}\label{gammaq}
\Gamma_\tau(x):=
\frac{q^{-\frac{x^2}{16}}}{\bigl(-q^{\frac{1}{2}(x+1)};q\bigr)_\infty}.
\end{equation}
Observe that 
$\Gamma_\tau(x)$ is a zero-free meromorphic function with
simple poles located at 
$-1+\tau^{-1}+2\mathbb{Z}_{\leq 0}+2\tau^{-1}\mathbb{Z}$.
It furthermore satisfies the difference equation
\begin{equation}\label{diffeq}
\Gamma_\tau(x+2)=2\cos(\pi(x+1)\tau/2)\Gamma_\tau(x).
\end{equation}
Note furthermore that for $\tau\in i\mathbb{R}_{>0}$, i.e. $0<q<1$, 
the function $\Gamma_\tau(z)$ satisfies
$\overline{\Gamma_\tau(x)}=\Gamma_\tau(\overline{x})$, where the bar stands for
the complex conjugate.
 
Observe that the above defined $q$-gamma type functions are not
$\tau^{-1}$-periodic. 
It is probably
for this reason that formulas in $q$-analysis are usually
expressed in terms of $q$-shifted factorials 
instead of $q$-gamma functions. For our present
purposes the expressions in terms of $q$-gamma type functions are
convenient because it clarifies the similarities
with the $|q|=1$ case.


\subsection{The gamma function for $|q|=1$.}
In this subsection we take $\tau\in\mathbb{R}_{<0}$, 
whence $q=q_\tau=\exp(2\pi i\tau)$
satisfies $|q|=1$. As in the previous subsection, we write
$q^u=\exp(2\pi i\tau u)$ for $u\in\mathbb{C}$. 

It is easy to verify that the integral
\begin{equation}
\gamma_\tau(z)=\frac{1}{2i}\int_0^{\infty}\frac{dy}{y}\left(\frac{z}{y}-
\frac{\sinh(\tau yz)}{\sinh(y)\sinh(\tau y)}\right)
\end{equation}
converges absolutely in the strip $|\hbox{Re}(z)|<1-\tau^{-1}$. 
For $z\in\mathbb{C}$ in this strip we set
\begin{equation}
G_\tau(z)=\exp(i\gamma_\tau(z)).
\end{equation}
Ruijsenaars' \cite[Sect. 3]{Ru1}
hyperbolic gamma function $G(z)=G(a_+,a_-;z)$ with $a_+,a_->0$
is related to $G_\tau$ by 
\[G_\tau(z)=G(a_+,a_-;ia_-z/2),\qquad \tau=-a_-/a_+.
\]
In the following proposition we recall 
some of Ruijsenaars' results \cite[Sect. 3]{Ru1}
on the hyperbolic gamma function.
\begin{prop}\label{Ruij}
{\bf (i)} The function $G_\tau(z)$ satisfies the difference
equation
\[G_\tau(z+2)=2\cos(\pi(z+1)\tau/2)G_\tau(z).
\]
In particular, $G_\tau(z)$ admits a meromorphic continuation to the complex
plane $\mathbb{C}$, which we again denote by $G_\tau(z)$.

{\bf (ii)} The zeros of $G_\tau(z)$ are located at
$1-\tau^{-1}+2\mathbb{Z}_{\geq 0}+2\tau^{-1}\mathbb{Z}_{\leq 0}$,
and the poles of $G_\tau(z)$ are located at
$-1+\tau^{-1}+2\mathbb{Z}_{\leq 0}+2\tau^{-1}\mathbb{Z}_{\geq 0}$. 

{\bf (iii)} $G_\tau(-z)G_\tau(z)=1$ and $G_\tau(z)=G_{\tau^{-1}}(-\tau z)$.

{\bf (iv)} The function $\gamma_\tau(z)$ has an analytic continuation to the
cut plane 
\[\mathbb{C}\setminus\{ (-\infty,-1+\tau^{-1}]\cup [1-\tau^{-1},\infty)\},
\]
which we again denote by $\gamma_\tau(z)$. 
Set $r=\max(1,-\tau^{-1})$ and choose $\epsilon>0$, then
\[|\mp \gamma_\tau(z)+\frac{\pi\tau z^2}{8}-\frac{\pi}{24}(\tau+\tau^{-1})|=
\mathcal{O}\bigl(\exp((\epsilon-\pi/r)|\hbox{Im}(z)|)\bigr), \qquad
\hbox{Im}(z)\rightarrow \pm\infty
\]
uniformly for $\hbox{Re}(z)$ in compacts of $\mathbb{R}$.
\end{prop}
{}From the explicit expression for 
$G_\tau(z)$ with $|\hbox{Re}(z)|<1-\tau^{-1}$
and the first order difference equation for $G_\tau$, we have
$\overline{G_\tau(z)}=G_\tau(\overline{z})$.

As Ruijsenaars verifies in 
\cite[Appendix A]{Ru2}, the hyperbolic gamma function is
a quotient of Barnes' double gamma function, and it
essentially coincides with Kurokawa's double sine function.

The hyperbolic gamma function is the important building block for
$q$-analysis with $|q|=1$. It was used in \cite{N} and \cite{NU}
to construct for $|q|=1$
explicit integral solutions of the $q$-Bessel difference
equation and of the $q$-hypergeometric difference equation. 
Ruijsenaars' \cite{Ru2} used hyperbolic
gamma functions to construct an eigenfunction for the Askey-Wilson second
order difference operator 
for $|q|=1$ as an explicit Barnes' type integral. In this
paper we construct an eigenfunction
of the Askey-Wilson second order difference operator for $|q|=1$ as an Euler
type integral, using representation theory of quantum groups.


\section{The Askey-Wilson second order difference operator and quantum groups}

Throughout this section we require that 
$\tau\in\mathbb{C}\setminus \frac{1}{2}\mathbb{Z}$
and we write $q=q_\tau=\exp(2\pi i\tau)$ 
for the corresponding deformation parameter.
The condition on $\tau$ implies $q\not=\pm 1$. 
As usual, we write $q^u=\exp(2\pi i\tau u)$ for $u\in\mathbb{C}$.
\begin{defi}
The quantum group $\mathcal{U}_q$ is the unital associative algebra
over $\mathbb{C}$ generated by $K^{\pm 1}$, $X^+$ and $X^-$, subject to
the relations
\begin{equation*}
\begin{split}
&KK^{-1}=K^{-1}K=1,\\
&KX^+=qX^+K,\qquad KX^-=q^{-1}X^-K,\\
&X^+X^--X^-X^+=\frac{K^2-K^{-2}}{q-q^{-1}}.
\end{split}
\end{equation*}
\end{defi}
It is well known that $\mathcal{U}_q$ has the structure of 
a Hopf-algebra, and as such it is a quantum deformation of the
universal enveloping algebra of the simple
Lie algebra $\mathfrak{s}\mathfrak{l}(2,\mathbb{C})$. The Hopf-algebra
structure does not play a significant role in the present paper, so
the definition of the Hopf-algebra structure is omitted here. I only want to
stipulate that the upcoming definition of Koornwinder's twisted
primitive element is motivated by its transformation behaviour 
under the action of the comultiplication of $\mathcal{U}_q$, see \cite{K}
for details. 

It is convenient to work with an extended version of $\mathcal{U}_q$, which
we define as follows.
Write $\mathcal{A}=\bigoplus_{x\in\mathbb{C}}
\mathbb{C}\,\widehat{x}$ for the group algebra of the additive 
group $(\mathbb{C},+)$.
Denote $\hbox{End}_{alg}(\mathcal{U}_q)$ for the 
unital algebra homomorphisms $\phi:
\mathcal{U}_q\rightarrow\mathcal{U}_q$. There exists an 
algebra homomorphism
\[\kappa: \mathcal{A}\rightarrow \hbox{End}_{alg}(\mathcal{U}_q),
\]
with $\kappa(\widehat{x})=\kappa_x\in\hbox{End}_{alg}(\mathcal{U}_q)$ 
for $x\in\mathbb{C}$
defined by
\[\kappa_x(K^{\pm 1})=K^{\pm 1},\qquad
\kappa_x(X^{\pm})=q^{\pm x}X^{\pm}.
\]
Note that for $m\in\mathbb{Z}$,
\[\kappa_{m}(X)=K^{m}XK^{-m}\qquad \forall X\in \mathcal{U}_q,
\]
so the au\-to\-mor\-phisms $\kappa_{x}$ generalize
the inner au\-to\-mor\-phisms $K^{m}(\,\cdot\,)K^{-m}$ of 
$\mathcal{U}_q$ ($m\in\mathbb{Z}$).
The ex\-ten\-ded al\-ge\-bra $\widetilde{\mathcal{U}}_q$ is now de\-fi\-ned as
follows.
\begin{defi}
The unital, associative algebra $\widetilde{\mathcal{U}}_q$ is the
vector space $\mathcal{A}\otimes \mathcal{U}_q$ with multiplication defined by
\[(\widehat{x}\otimes X)(\widehat{y}\otimes Y)=\widehat{(x+y)}\otimes
\kappa_{-y}(X)Y,\qquad \forall\, x,y\in\mathbb{C},\qquad
\forall\,X,Y\in\mathcal{U}_q.
\]
The unit element is $\widehat{0}\otimes 1$.
\end{defi}
Observe that $\mathcal{A}$ and $\mathcal{U}_q$ embed
as algebras in $\widetilde{\mathcal{U}}_q$ by the formulas
\[a\mapsto a\otimes 1,\qquad X\mapsto \widehat{0}\otimes X
\]
for $a\in\mathcal{A}$ and $X\in\mathcal{U}_q$. We will use these canonical 
embeddings to identify the algebras $\mathcal{U}_q$ and $\mathcal{A}$
with their images in $\widetilde{\mathcal{U}}_q$. The commutation
relations between $\mathcal{A}$ and $\mathcal{U}_q$ within
$\widetilde{\mathcal{U}}_q$ then become
\[\widehat{x}\,X=\widehat{x}\otimes X=\kappa_x(X)\,\widehat{x},
\qquad 
x\in\mathbb{C},\,\, X\in\mathcal{U}_q.
\]
The {\it quantum Casimir element}, defined by 
\begin{equation}\label{Omega}
\Omega:=X^+X^-+\frac{q^{-1}K^2+qK^{-2}-2}{(q-q^{-1})^2}\in\mathcal{U}_q,
\end{equation}
is an algebraic generator of the center
$\mathcal{Z}(\mathcal{U}_q)$ of $\mathcal{U}_q$. Note that
$\Omega$ is also in the center 
$\mathcal{Z}(\widetilde{\mathcal{U}}_q)$ of the extended algebra
$\widetilde{\mathcal{U}}_q$. 

We now consider the following explicit
realization of $\widetilde{\mathcal{U}}_q$.
Let $\mathcal{M}$ be the space 
of meromorphic functions on the complex plane $\mathbb{C}$.
For any $\lambda\in\mathbb{C}$, the assignment
\begin{equation*}
\begin{split}
\bigl(\pi_\lambda(X^{\pm})f\bigr)(z)&=
q^{\pm z}\left(\frac{q^{-\frac{1}{2}-i\lambda}f(z\mp 1)
-q^{\frac{1}{2}+i\lambda}f(z\pm 1)}{q^{-1}-q}\right),\\
\bigl(\pi_\lambda(K^{\pm 1})f\bigr)(z)&=f(z\pm 1),\\
\bigl(\pi_\lambda(\widehat{x})f\bigr)(z)&=f(z+x),\qquad x\in\mathbb{C},\\
\end{split}
\end{equation*}
uniquely extends to a representation of $\widetilde{\mathcal{U}}_q$
on $\mathcal{M}$. The quantum Casimir element $\Omega$ acts as
\begin{equation}\label{Omegaaction}
\pi_\lambda(\Omega)=
\left(\frac{q^{i\lambda}-q^{-i\lambda}}{q-q^{-1}}\right)^2\,
\hbox{Id}.
\end{equation}
Observe furthermore that $\pi_\lambda(\widehat{m})=
\pi_\lambda(K^{m})$ for all $m\in\mathbb{Z}$.

In the present paper the quantum group input to the theory of $q$-special
functions is based
on the explicit connection 
between the radial part of the quantum Casimir element
$\Omega$ and the second order Askey-Wilson difference operator.  
Here the {\it Askey-Wilson second order difference operator} $\mathcal{D}=
\mathcal{D}^{a,b,c,d}$, depending on four
parameters $(a,b,c,d)$ called the Askey-Wilson parameters, is
defined by
\begin{equation}\label{D}
(\mathcal{D}f)(x)=A(x)(f(x+2)-f(x))+A(-x)(f(x-2)-f(x)),
\end{equation}
with $A(x)=A(x;a,b,c,d)$ the explicit function
\[A(x)=\frac{(1-q^{a+x})(1-q^{b+x})(1-q^{c+x})(1-q^{d+x})}
{(1-q^{2x})(1-q^{2+2x})},
\]
cf. \cite{AW}.
The radial part of $\Omega$ is computed with respect to
elements $Y_{\rho}-\mu_\alpha(\rho)\,1\in\mathcal{U}_q$ for
$\alpha,\rho\in\mathbb{C}$, where
$Y_{\rho}$ is Koornwinder's \cite{K} {\it twisted primitive element},
\begin{equation}\label{Y}
Y_{\rho}=q^{\frac{1}{2}}X^+K-q^{-\frac{1}{2}}X^-K+
\left(\frac{q^{-\rho}+q^{\rho}}{q^{-1}-q}\right)(K^2-1)
\end{equation}
and $\mu_\alpha(\rho)$ is the constant
\begin{equation}\label{mu}
\mu_\alpha(\rho)=
\left(\frac{q^\rho(1-q^\alpha)+q^{-\rho}(1-q^{-\alpha})}
{q-q^{-1}}\right).
\end{equation}
Consider the five dimensional space
\begin{equation}
\mathcal{U}_q^1=
\hbox{span}_{\mathbb{C}}\{X^+,X^-,K,K^{-1},1\}\subset \mathcal{U}_q.
\end{equation}
The radial part computation of $\Omega$ leads to the following result.
\begin{prop}\label{radial}
Let $\rho,\sigma,\alpha,\beta\in\mathbb{C}$.
For all $x\in\mathbb{C}$,
\[\widehat{x}\,\Omega K=\widehat{x}\,\Omega(x)K\qquad \hbox{mod}\,\,\,\,
(Y_{\rho}-\mu_\alpha(\rho))\,\widehat{x}\,\mathcal{U}_q^1+
\widehat{x}\,\mathcal{U}_q^1\,(Y_{\sigma}-\mu_\beta(\sigma))
\]
with $\Omega(x)=\Omega(x;\alpha,\rho,\beta,\sigma)$ given explicitly by
\[\Omega(x)=\frac{q^{\beta-1}}{(q-q^{-1})^2}
\left\{B(x)K^2+\bigl(C(x)+(1-q^{1-\beta})^2\bigr)1+D(x)K^{-2}\right\}
\]
where 
\begin{equation*}
\begin{split}
B(x)&=B(x;\alpha,\rho,\beta,\sigma)=
q^{-\beta}\frac{(1-q^{a+x})(1-q^{2-a+x})(1-q^{b+x})(1-q^{2-b+x})}
{(1-q^{2x})(1-q^{2+2x})},\\
C(x)&=C(x;\alpha,\rho,\beta,\sigma)=-A(x;a,b,c,d)-A(-x;a,b,c,d),\\
D(x)&=D(x;\alpha,\rho,\beta,\sigma)=
q^{-\beta}\frac{(1-q^{c-x})(1-q^{2-c-x})(1-q^{d-x})(1-q^{2-d-x})}
{(1-q^{-2x})(1-q^{2-2x})}.
\end{split}
\end{equation*}
Here
the Askey-Wilson parameters $(a,b,c,d)$ are related to the parameters
$\alpha,\beta,\rho,\sigma$ by 
\begin{equation}\label{parlinks}
(a,b,c,d)=(1+\rho+\sigma,1-\rho+\sigma,
1+\alpha+\rho-\beta-\sigma,1-\alpha-\rho-\beta-\sigma).
\end{equation}
\end{prop}
\begin{proof}
The proof generalizes the radial part computation by
Koornwinder \cite{K}, where the case $\alpha=\beta=0$
and $x\in\mathbb{Z}$ is considered (see also 
Noumi \& Mimachi \cite{NM} and Koelink \cite{Ko} 
for extensions to discrete values of $\alpha$ and $\beta$). 
In the present set-up the computation is a bit more complex,
and we gather more precise information on the remainder. 
For the convenience of the reader,
I have included the main steps of the proof as appendix.
\end{proof}
Unless specified otherwise, we assume that the Askey-Wilson parameters
$(a,b,c,d)$ are related to the four parameters $(\alpha,\rho,\beta,\sigma)$
by \eqref{parlinks}.

Proposition \ref{radial} allows us to identify specific eigenfunctions of 
$\pi_\lambda(\Omega)$ with eigenfunctions of the second order difference
operator
\begin{equation}
(\mathcal{L}f)(x)=(\mathcal{L}^{\alpha,\rho,\beta,\sigma}f)(x)=
B(x)f(x+2)+C(x)f(x)+D(x)f(x-2).
\end{equation}
 The second order difference operator $\mathcal{L}$
is gauge equivalent to the Askey-Wilson second order difference equation
$\mathcal{D}=\mathcal{D}^{a,b,c,d}$, since
$\mathcal{L}=\Delta\circ \mathcal{D}\circ \Delta^{-1}$ with
$\Delta(x)=\Delta(x;a,b,c,d)$ any meromorphic function 
satisfying the difference
equation
\begin{equation}\label{diffDelta}
\begin{split}
\Delta(x+2)&=\frac{\sin(\pi(c+x)\tau)\sin(\pi(d+x)\tau)}
{\sin(\pi(2-a+x)\tau)\sin(\pi(2-b+x)\tau)}\,\Delta(x)\\
&=\frac{(1-q^{c+x})(1-q^{d+x})}
{(1-q^{2-a+x})(1-q^{2-b+x})}\,q^\beta\,\Delta(x).
\end{split}
\end{equation}

To make use of the above radial part computation, we first need to 
construct explicit
eigenfunctions of $\pi_\lambda(Y_{\rho})$ with eigenvalue $\mu_{\alpha}(\rho)$.
As we will see in the proof of the following proposition,
the operator $\pi_\lambda(Y_{\rho})$ is a  
first order difference operator, hence eigenfunctions of 
$\pi_\lambda(Y_{\rho})$ admit the following simple characterization.
\begin{prop}\label{diffone}
Let $\alpha,\rho\in\mathbb{C}$.
A meromorphic function $f\in\mathcal{M}$ is an eigenfunction of
$\pi_\lambda(Y_{\rho})$ with eigenvalue $\mu_\alpha(\rho)$
if and only if 
\[f(z+2)=\frac{\sin(\pi(-i\lambda-\alpha-\rho+z)\tau)
\sin(\pi(-i\lambda+\alpha+\rho+z)\tau)}
{\sin(\pi(i\lambda-\rho+1+z)\tau)
\sin(\pi(i\lambda+\rho+1+z)\tau)}\,f(z).
\]
\end{prop}
\begin{proof}
A direct computation shows that
$\pi_\lambda(Y_{\rho})\in\hbox{End}(\mathcal{M})$ is the explicit first
order difference operator
\begin{equation*}
\begin{split}
\bigl(\pi_\lambda(Y_{\rho})f\bigr)(z)=
\frac{q^{-\rho}}{q-q^{-1}}
&\left\{(q^{\rho-1-i\lambda-z}-1)(1-q^{\rho+1+i\lambda+z})f(z+2)\right.\\
&\qquad\qquad\left.-(q^{\rho+i\lambda-z}-1)(1-
  q^{\rho-i\lambda+z})f(z)\right\},
\end{split}
\end{equation*}
and, more generally,
\begin{equation*}
\begin{split}
\bigl((\pi_\lambda(Y_{\rho})-\mu_\alpha(\rho))f\bigr)(z)=
\frac{q^{-\rho}}{q-q^{-1}}
&\left\{(q^{\rho-1-i\lambda-z}-1)(1-
  q^{\rho+1+i\lambda+z})f(z+2)\right.\\
&\quad\qquad\left.-(q^{\alpha+\rho+i\lambda-z}-1)(1-
  q^{\alpha+\rho-i\lambda+z})q^{-\alpha}f(z)\right\}.
\end{split}
\end{equation*}
The eigenvalue equation $\pi_\lambda(Y_{\rho})f=
\mu_\alpha(\rho)f$ is thus equivalent to the first order difference equation
\[ f(z+2)=\frac{(1-q^{\alpha+\rho+i\lambda-z})
(1-q^{\alpha+\rho-i\lambda+z})}{(1-q^{\rho-1-i\lambda-z})
(1-q^{\rho+1+i\lambda+z})}q^{-\alpha}f(z).
\]
Rewriting this formula yields the desired result.
\end{proof}


\section{The Askey-Wilson function for $0<q<1$.}
In this section we take $\tau\in i\mathbb{R}_{>0}$, so that
$0<q=q_\tau=\exp(2\pi i\tau)<1$. The assignment
\begin{equation}
(K^{\pm 1})^*=K^{\pm 1},\qquad (X^{\pm })^*=-X^{\mp}
\end{equation}
uniquely extends to a unital, anti-linear, anti-algebra involution
on $\mathcal{U}_q$. This particular choice of $*$-structure 
corresponds classically to 
choosing the real form $\mathfrak{s}\mathfrak{u}(1,1)$ of
$\mathfrak{s}\mathfrak{l}(2,\mathbb{C})$. 

The elements $K^m$ ($m\in\mathbb{Z}$), the quantum Casimir element 
$\Omega$ (see \eqref{Omega}), 
and the special family $Y_\rho$ ($\rho\in\mathbb{R}$)
of Koornwinder's twisted primitive elements \eqref{Y} are $*$-selfadjoint
elements in $\mathcal{U}_q$. The eigenvalue $\mu_\alpha(\rho)$ (see \eqref{mu})
is real for $\alpha,\rho\in\mathbb{R}$.
We consider now an explicit 
$*$-unitary pairing for the representation $\pi_\lambda$.
\begin{lem}\label{selfadjoint1}
Let $\lambda\in\mathbb{R}$.
Suppose that $f,g\in\mathcal{M}$ are $\tau^{-1}$-periodic and 
analytic on the strip $\{z\in\mathbb{C}\, 
|\, |\hbox{Re}(z)|\leq 1\}$. Then
\[\langle \pi_\lambda(X)f,g\rangle=\langle f,\pi_\lambda(X^*)g\rangle,\qquad
\forall X\in \mathcal{U}_q^1,
\]
with the pairing $\langle \cdot,\cdot\rangle$ defined by
\[\langle f,g\rangle:=\int_0^1f(y/\tau)\overline{g(y/\tau)}dy.
\]
\end{lem}
\begin{proof}
This is an easy verification for the basis elements $1,K^{\pm 1}, X^{\pm}$
of $\mathcal{U}_q^1$, using Cauchy's Theorem to shift contours.
\end{proof}
\begin{rema}
This lemma can be applied recursively. 
Let $f,g\in\mathcal{M}$ be $\tau^{-1}$-periodic and
analytic on the strip $\{z\in\mathbb{C}\,|\,|\hbox{Re}(z)|\leq k\}$ with 
$k\in\mathbb{Z}_{>0}$.
For any $X=X_1X_2\cdots X_m\in \mathcal{U}_q$ with $m\leq k$ and 
$X_i\in \mathcal{U}_q^1$,
\[\langle \pi_\lambda(X)f,g\rangle=\langle f,\pi_\lambda(X^*)g\rangle.
\]
In particular, the subspace 
$\mathcal{O}_{\tau^{-1}}$ of entire, $\tau^{-1}$-periodic
functions is an $*$-unitary
$\pi_\lambda$-invariant subspace of $\mathcal{M}$ with respect to the pairing
$\langle \cdot,\cdot\rangle$. 
This subspace is the algebraic version of the 
principal series representation of the $*$-algebra $(\mathcal{U}_q,*)$.
\end{rema}
To combine this lemma with the radial part computation of
the quantum Casimir element (see Proposition \ref{radial}), 
we need to construct meromorphic $\tau^{-1}$-periodic
eigenfunctions of $\pi_\lambda(Y_\rho)$ ($\rho\in\mathbb{R}$) which are
analytic in a large enough strip around the imaginary axis. We claim that
the meromorphic function
\begin{equation}\label{f}
f_\lambda(z)=f_\lambda(z;\alpha,\rho):=
\frac{\Gamma_{2\tau}\bigl(-1-\frac{1}{2\tau}+\alpha+\rho-i\lambda+z\bigr)
\Gamma_{2\tau}\bigl(-\frac{1}{2\tau}+\rho-i\lambda-z\bigr)}
{\Gamma_{2\tau}\bigl(1-\frac{1}{2\tau}+\alpha+\rho+i\lambda-z\bigr)
\Gamma_{2\tau}\bigl(-\frac{1}{2\tau}+\rho+i\lambda+z\bigr)}
\end{equation}
meets these criteria for special values of the parameters. 
By Proposition \ref{diffone} and the difference equation for 
$\Gamma_\tau$ it follows that $f_\lambda(z;\alpha,\rho)$
is an eigenfunction of $\pi_\lambda(Y_\rho)$ with eigenvalue
$\mu_\alpha(\rho)$ for $\alpha,\rho,\lambda\in\mathbb{C}$.
Writing $f_\lambda(z;\alpha,\rho)$ 
in terms of $q$-shifted factorials leads to the
expression
\begin{equation}
f_\lambda(z;\alpha,\rho)=C\,\frac{\bigl(q^{2+\alpha+\rho+i\lambda-z},
q^{1+\rho+i\lambda+z};q^2\bigr)_{\infty}}
{\bigl(q^{\alpha+\rho-i\lambda+z},q^{1+\rho-i\lambda-z};q^2\bigr)_{\infty}}\,
q^{-\frac{\alpha z}{2}}
\end{equation}
for some nonzero constant $C$ (independent of $z$), hence $f_\lambda(z)$
is $\tau^{-1}$-periodic if and only if 
$\alpha\in 2\mathbb{Z}$. Furthermore, observe
that the poles of $f_\lambda(z)$ are located at
\begin{equation}\label{polesg}
i\lambda-\alpha-\rho+2\mathbb{Z}_{\leq 0}+\mathbb{Z}\tau^{-1},\qquad
-i\lambda+\rho+1+2\mathbb{Z}_{\geq 0}+\mathbb{Z}\tau^{-1},
\end{equation}
so $f_\lambda(z)$ is analytic on the strip
$\{z\in\mathbb{C}\,|\,|\hbox{Re}(z)|\leq\rho\}$ 
when $\alpha>0$, $\rho\geq 0$ and $\lambda\in\mathbb{R}$.
\begin{defi}\label{defin1}
Let $\lambda\in\mathbb{R}$, $\alpha,\beta\in 2\mathbb{Z}_{>0}$ and
$\rho,\sigma\in\mathbb{R}_{\geq 3}$. For $x\in\mathbb{C}$ with
$|\hbox{Re}(x)|\leq 2$ we define 
$\phi_\lambda(x)=\phi_\lambda(x;\alpha,\rho,\beta,\sigma)$
by
\[\phi_\lambda(x):=
\langle \pi_\lambda(\widehat{x}K)f_\lambda(\cdot;\beta,\sigma),
f_\lambda(\cdot;\alpha,\rho)\rangle.
\]
\end{defi}
Note that the matrix coefficient $\phi_\lambda(x)$ is given explicitly by 
\begin{equation}\label{phi1}
\begin{split}
\phi_\lambda(x)=
\int_0^1&\frac{\Gamma_{2\tau}\bigl(-\frac{1}{2\tau}+\beta+\sigma-i\lambda+x+
\frac{y}{\tau}\bigr)
\Gamma_{2\tau}\bigl(-1-\frac{1}{2\tau}+\sigma-i\lambda-x-\frac{y}{\tau}\bigr)}
{\Gamma_{2\tau}\bigl(-\frac{1}{2\tau}+\beta+
\sigma+i\lambda-x-\frac{y}{\tau}\bigr)
\Gamma_{2\tau}\bigl(1-\frac{1}{2\tau}+
\sigma+i\lambda+x+\frac{y}{\tau}\bigr)}\\
&\quad\qquad\times
\frac{\Gamma_{2\tau}\bigl(-1+\frac{1}{2\tau}+
\alpha+\rho+i\lambda-\frac{y}{\tau}\bigr)
\Gamma_{2\tau}\bigl(\frac{1}{2\tau}+\rho+i\lambda+\frac{y}{\tau}\bigr)}
{\Gamma_{2\tau}\bigl(1+\frac{1}{2\tau}+
\alpha+\rho-i\lambda+\frac{y}{\tau}\bigr)
\Gamma_{2\tau}\bigl(\frac{1}{2\tau}+\rho-i\lambda-\frac{y}{\tau}\bigr)}\,dy,
\end{split}
\end{equation}
and that $\phi_\lambda(x)$ is analytic on the strip
$\{x\in\mathbb{C}\,|\,|\hbox{Re}(x)|\leq 2\}$.
The quantum group interpretation of this explicit integral leads to the
following result.
\begin{thm}\label{main1}
Let $\lambda\in\mathbb{R}$, $\alpha,\beta\in 2\mathbb{Z}_{>0}$
and $\rho,\sigma\in\mathbb{R}_{\geq 3}$.
The matrix coefficient $\phi_\lambda(x)=
\phi_\lambda(x;\alpha,\rho,\beta,\sigma)$ 
satisfies the second order difference equation
\[\bigl(\mathcal{L}\phi_\lambda\bigr)(x)=E(\lambda)\phi_\lambda(x)
\]
for generic $x\in i\mathbb{R}$,
with the eigenvalue $E(\lambda)=E(\lambda;\beta)$ given by
\begin{equation}\label{E}
E(\lambda)=-1-q^{2-2\beta}+q^{1-\beta}(q^{2i\lambda}+q^{-2i\lambda}).
\end{equation}
\end{thm}
\begin{proof}
For the duration of the proof we use the shorthand notations
$f(z)=f_\lambda(z;\beta,\sigma)$ and $g(z)=f_\lambda(z;\alpha,\rho)$.
By the conditions on the parameters, 
$Y_\rho$ is $*$-selfadjoint, $\mu_\alpha(\rho)$ is real
and the meromorphic functions $f(z)$ and $g(z)$ are analytic on the strip
$\{z\in\mathbb{C}\,|\,|\hbox{Re}(z)|\leq 3\}$. 
By Lemma \ref{selfadjoint1} we thus 
obtain for any $X\in \mathcal{U}_q^1$ and $x\in i\mathbb{R}$,
\[\langle \pi_\lambda((Y_\rho-\mu_\alpha(\rho))\widehat{x}X)f,g\rangle=
\langle \pi_\lambda(\widehat{x}X)f,
\pi_\lambda(Y_\rho-\mu_\alpha(\rho))g\rangle=0,
\]
and obviously also
$\langle \pi_\lambda(\widehat{x}X(Y_\sigma-\mu_\beta(\sigma)))f,g\rangle=0$.
We conclude from Proposition \ref{radial} that
\begin{equation*}
\begin{split}
\langle \pi_\lambda(\widehat{x}\,\Omega K)f,g\rangle
&=\langle \pi_\lambda(\widehat{x}\,\Omega(x)Kf,g\rangle\\
&=\frac{q^{\beta-1}}{(q-q^{-1})^2}\,\bigl(\mathcal{L}\phi_\lambda\bigr)(x)
+\frac{q^{\beta-1}(1-q^{1-\beta})^2}{(q-q^{-1})^2}\,\phi_\lambda(x).
\end{split}
\end{equation*}
On the other hand, \eqref{Omegaaction} implies that
\[\langle \pi_\lambda(\widehat{x}\,\Omega K)f,g\rangle=
\left(\frac{q^{i\lambda}-q^{-i\lambda}}{q-q^{-1}}\right)^2\phi_\lambda(x),
\]
hence $\bigl(\mathcal{L}\phi_\lambda\bigr)(x)=E(\lambda)\phi_\lambda(x)$.
\end{proof}
For the comparison with the results for $|q|=1$, see Section 5, 
it is convenient to note that the matrix coefficient $\phi_\lambda(x)$
can be rewritten as
\begin{equation}\label{phi2}
\begin{split}
\phi_\lambda(x)=C\,
\int_0^1&\frac{\Gamma_{2\tau}\bigl(\frac{1}{2\tau}+\beta+\sigma-i\lambda+x+
\frac{y}{\tau}\bigr)
\Gamma_{2\tau}\bigl(-1+\frac{1}{2\tau}+\sigma-i\lambda-x-\frac{y}{\tau}\bigr)}
{\Gamma_{2\tau}\bigl(-\frac{1}{2\tau}+\beta+\sigma+i\lambda-x-
\frac{y}{\tau}\bigr)
\Gamma_{2\tau}\bigl(1-\frac{1}{2\tau}+\sigma+i\lambda+x+\frac{y}{\tau}\bigr)}\\
&\quad\qquad\times
\frac{\Gamma_{2\tau}\bigl(-1-\frac{1}{2\tau}+
\alpha+\rho+i\lambda-\frac{y}{\tau}\bigr)
\Gamma_{2\tau}\bigl(-\frac{1}{2\tau}+\rho+i\lambda+\frac{y}{\tau}\bigr)}
{\Gamma_{2\tau}\bigl(1-\frac{1}{2\tau}+
\alpha+\rho-i\lambda+\frac{y}{\tau}\bigr)
\Gamma_{2\tau}\bigl(-\frac{1}{2\tau}+\rho-i\lambda-\frac{y}{\tau}\bigr)}\,dy
\end{split}
\end{equation}
for some $x$-independent nonzero 
constant $C$. Formula \eqref{phi2} follows from
\eqref{phi1} and the difference equation
\[\Gamma_{2\tau}\bigl(x+\frac{1}{2\tau}\bigr)=
\exp\bigl(-\frac{\pi ix}{2}\,\bigr)\Gamma_{2\tau}\bigl(x-\frac{1}{2\tau}\bigr).
\]
\begin{cor}\label{D1}
Let $\lambda\in\mathbb{R}$, $\alpha,\beta\in 2\mathbb{Z}_{>0}$
and $\rho,\sigma\in\mathbb{R}_{\geq 3}$.
Let the Askey-Wilson parameters $(a,b,c,d)$ be given by
\eqref{parlinks}.
The function $F_\lambda(x)=F_\lambda(x;a,b,c,d)$ defined by
$F_\lambda(x)=\Delta(x)^{-1}\phi_\lambda(x)$, 
with $\Delta(x)=\Delta(x;a,b,c,d)$
the $\tau^{-1}$-periodic, meromorphic function
\begin{equation}\label{Delta}
\Delta(x)=\frac{\Gamma_{2\tau}\bigl(-1+\frac{1}{2\tau}+a-x\bigr)
\Gamma_{2\tau}\bigl(-1+\frac{1}{2\tau}+b-x\bigr)}
{\Gamma_{2\tau}\bigl(1+\frac{1}{2\tau}-c-x\bigr)
\Gamma_{2\tau}\bigl(1+\frac{1}{2\tau}-d-x\bigr)},
\end{equation}
satisfies the Askey-Wilson second order difference equation
\[(\mathcal{D}F_\lambda)(x)=E(\lambda)F_\lambda(x)
\]
for generic $x\in i\mathbb{R}$ and generic $\rho$ and $\sigma$. 
\end{cor}
\begin{proof}
Note that $\Delta(x)$ can be rewritten as
\begin{equation}\label{Delta2}
\Delta(x)=\frac{\bigl(q^{2-c-x},q^{2-d-x};q^2\bigr)_{\infty}}
{\bigl(q^{a-x},q^{b-x};q^2\bigr)_{\infty}}\,q^{-\frac{\beta x}{2}}
\end{equation}
for some nonzero ($x$-independent) constant $C$.
Since $\beta\in 2\mathbb{Z}_{>0}$, the gauge factor $\Delta(x)$ is 
$\tau^{-1}$-invariant.
Furthermore, $\Delta(x)^{-1}$
is regular at $x\in\pm 2+i\mathbb{R}$ and $x\in i\mathbb{R}$
under generic conditions on the parameters $\rho$ and $\sigma$.
The proof is completed by observing that $\Delta(x)$ satisfies the difference
equation \eqref{diffDelta}. 
\end{proof}
We end this section by extending these results to continuous 
parameters $\alpha$ and $\beta$. 
Changing integration variable
and substituting the expression for $\Gamma_{2\tau}$ in terms of $q$-shifted
factorials, we can rewrite 
the matrix coefficient $\phi_\lambda(x)$ (see \eqref{phi1}) as
\begin{equation}\label{intformsu110}
\phi_\lambda(x)=C\,\frac{q^{-\frac{\beta x}{2}}}{2\pi i}
\int_\mathbb{T} 
\frac{\bigl(q^{1+\beta+\sigma+i\lambda-x}/z,q^{2+\sigma+i\lambda+x}z,
q^{2+\alpha+\rho-i\lambda}z, q^{1+\rho-i\lambda}/z;q^2\bigr)_{\infty}}
{\bigl(q^{\sigma-i\lambda-x}/z,q^{1+\beta+\sigma-i\lambda+x}z,
  q^{1+\rho+i\lambda}z, q^{\alpha+\rho+i\lambda}/z;q^2\bigr)_{\infty}}
z^{\frac{\alpha-\beta}{2}}\frac{dz}{z}
\end{equation}
for some $x$-independent nonzero constant $C$, where 
$\mathbb{T}=\{z\in \mathbb{C} \, | \, |z|=1\}$ is the positively oriented
unit circle in the complex plane. 
To allow $\alpha,\beta$ to be continuous parameters,
we need to get rid of the term 
$z^{\frac{\alpha-\beta}{2}}$ in the integrand. This can be
achieved by rewriting $\phi_\lambda(x)$ as
\begin{equation}\label{intformsu11}
\begin{split}
\phi_\lambda(x)=C\,
\frac{q^{-\frac{\beta x}{2}}}{2\pi i}\int_{\mathbb{T}}
&\frac{\bigl(q^{1+\beta+\sigma+i\lambda-x}/z, 
q^{1+\alpha+\rho-\beta-i\lambda}/z;q^2\bigr)_{\infty}}
{\bigl(q^{1-\rho+i\lambda}z,q^{1+\rho+i\lambda}z;q^2\bigr)_{\infty}}\\
&\qquad\times\frac{\bigl(q^{1-\alpha-\rho+
\beta+i\lambda}z,q^{2+\sigma+i\lambda+x}z, 
q^{2+\alpha+\rho-i\lambda}z;q^2\bigr)_\infty}
{\bigl(q^{\alpha+\rho+i\lambda}/z,
q^{\sigma-i\lambda-x}/z,q^{1+\beta+\sigma-i\lambda+x}z
;q^2\bigr)_{\infty}}\frac{dz}{z},
\end{split}
\end{equation}
with $C$ again some (different) irrelevant $x$-independent nonzero constant.
The integral formula \eqref{intformsu11} follows from
\eqref{intformsu110} by substitution of the identity
\[z^{\frac{\alpha-\beta}{2}}\bigl(q^{1+\rho-i\lambda}/z;q^2\bigr)_{\infty}=
(-q^{i\lambda-\rho})^{\frac{\beta-\alpha}{2}}q^{\frac{(\beta-\alpha)^2}{4}}\,
\frac{\bigl(q^{1-\alpha-\rho+\beta+i\lambda}z,
q^{1+\alpha+\rho-\beta-i\lambda}/z;q^2\bigr)_{\infty}}
{\bigl(q^{1-\rho+i\lambda}z;q^2\bigr)_{\infty}},
\]
which in turn is a direct consequence of the functional equation
\[\theta(q^{2k}z)=(-z)^{-k}q^{-k^2}\theta(z),\qquad k\in\mathbb{Z}
\]
for the modified Jacobi theta function
$\theta(z)=\bigl(qz,q/z;q^2\bigr)_{\infty}$.
By \eqref{Delta2}, the eigenfunction $F_\lambda(x)=F_\lambda(x;a,b,c,d)$ 
of the Askey-Wilson
second order difference equation $\mathcal{D}$ (see Corollary \ref{D1}) 
is equal to
\begin{equation}\label{F}
\begin{split}
\mathcal{F}_\lambda(x)=\mathcal{F}_\lambda(x;a,b,c,d)&:=
\frac{\bigl(q^{1+\rho+\sigma-x},q^{1-\rho+\sigma-x};q^2\bigr)_{\infty}}
{\bigl(q^{1-\alpha-\rho+\beta+\sigma-x},
q^{1+\alpha+\rho+\beta+\sigma-x};q^2\bigr)_{\infty}}\\
&\,\,\,\times \frac{1}{2\pi i}\int_{\mathbb{T}}
\frac{\bigl(q^{1+\beta+\sigma+i\lambda-x}/z, 
q^{1+\alpha+\rho-\beta-i\lambda}/z;q^2\bigr)_{\infty}}
{\bigl(q^{1-\rho+i\lambda}z,q^{1+\rho+i\lambda}z;q^2\bigr)_{\infty}}\\
&\qquad\qquad\times\frac{\bigl(q^{1-\alpha-\rho+
\beta+i\lambda}z,q^{2+\sigma+i\lambda+x}z, 
q^{2+\alpha+\rho-i\lambda}z;q^2\bigr)_\infty}
{\bigl(q^{\alpha+\rho+i\lambda}/z,
q^{\sigma-i\lambda-x}/z,q^{1+\beta+\sigma-i\lambda+x}z
;q^2\bigr)_{\infty}}\frac{dz}{z}
\end{split}
\end{equation}
up to some nonzero $x$-independent multiplicative constant. 

Define for the Askey-Wilson parameters $(a,b,c,d)$ given by
\eqref{parlinks}, dual Askey-Wilson
parameters $(\widetilde{a},\widetilde{b},\widetilde{c},\widetilde{d})$
by
\[(\widetilde{a},\widetilde{b},\widetilde{c},\widetilde{d})=
(1-\beta,1+\beta+2\sigma,1+\alpha+2\rho,1-\alpha).
\]
This notion of dual Askey-Wilson parameters coincides with the
notion of dual parameters as used in e.g. \cite{KS2}.

In the following theorem we express 
$\mathcal{F}_\lambda(x)$ in terms of basic hypergeometric series by shrinking
the radius of the integration circle 
$\mathbb{T}$ to zero while picking up residues.
Recall that the very-well-poised ${}_8\phi_7$ series is defined by
\[{}_8W_7\bigl(u;b_1,b_2,b_3,b_4,b_5;q,z\bigr)=
\sum_{k=0}^{\infty}
\frac{(1-uq^{2k})\bigl(u,b_1,b_2,b_3,b_4,b_5;q\bigr)_kz^k}
{(1-u)\bigl(q,qu/b_1,qu/b_2,qu/b_3,qu/b_4,qu/b_5;q\bigr)_k}
\]
for $|z|<1$, see \cite{GR}. 
\begin{thm}
Let $\lambda\in \mathbb{R}$, $\alpha,\beta\in\mathbb{R}_{>0}$ and
$\rho,\sigma\in\mathbb{R}_{\geq 3}$. 
Under the parameter correspondence \eqref{parlinks},
the function $\mathcal{F}_\lambda(x)=\mathcal{F}_\lambda(x;a,b,c,d)$ 
given by \eqref{F} can be expressed
in terms of basic hypergeometric series as
\begin{equation}\label{hs}
\begin{split}
\mathcal{F}_\lambda(x)=C\,&\frac{\bigl(q^{2+a-\widetilde{d}+2i\lambda+x},
q^{2+a-\widetilde{d}+2i\lambda-x};q^2\bigr)_{\infty}}
{\bigl(q^{2-d+x},q^{2-d-x};q^2\bigr)_{\infty}}\\
&\times {}_8W_7\bigl(
q^{-2+\widetilde{a}+\widetilde{b}+\widetilde{c}+2i\lambda};
q^{a+x},q^{a-x},q^{\widetilde{a}+2i\lambda},q^{\widetilde{b}+2i\lambda},
q^{\widetilde{c}+2i\lambda};q^2,q^{2-\widetilde{d}-2i\lambda}\bigr)
\end{split}
\end{equation}
with the \textup{(}irrelevant\textup{)} 
generically nonzero, $x$-independent constant $C$
given by
\[C=\frac{\bigl(q^{2+2\sigma},q^{2+\alpha+2\rho-\beta},
q^{2+\alpha+2\rho+\beta+2\sigma}, q^{1+\beta+2i\lambda},
q^{1+\alpha-2i\lambda};q^2\bigr)_{\infty}}
{\bigl(q^2,q^{1+\alpha+2i\lambda}, q^{1+\alpha+2\rho+2i\lambda},
q^{1+\beta+2\sigma-2i\lambda},q^{3+\alpha+2\rho+2\sigma+2i\lambda};
q^2\bigr)_{\infty}}.
\]
Furthermore, $\mathcal{F}_\lambda(x)$ is $\tau^{-1}$-periodic and 
satisfies the Askey-Wilson difference equation
$\bigl(\mathcal{D}\mathcal{F}_\lambda\bigr)(x)=
E(\lambda)\mathcal{F}_\lambda(x)$ 
for generic $x\in i\mathbb{R}$.
\end{thm}
\begin{proof}
The expression for $\mathcal{F}_\lambda$ 
follows by shrinking the radius of the integration contour
$\mathbb{T}$ to zero while picking up residues. 
It is actually a special case of 
\cite[Exerc. 4.4, p.122]{GR}, 
in which one should replace the base $q$ by $q^2$ and the parameters
$(a,b,c,d,f,g,h,k)$  by
\[\bigl(q^{1-\beta-\sigma-i\lambda+x},q^{1-\rho+i\lambda},
q^{\alpha+\rho+i\lambda},q^{\sigma-i\lambda-x},q^{2+\sigma+i\lambda+x},
q^{1+\rho+i\lambda},q^{1+\beta+\sigma-i\lambda+x},
q^{1-\beta-2i\lambda}\bigr).
\]
We have seen that the difference equation
$\bigl(\mathcal{D}\mathcal{F}_\lambda\bigr)(x)=
E(\lambda)\mathcal{F}_\lambda(x)$ for generic $x\in i\mathbb{R}$
is valid under the extra assumption 
$\alpha,\beta\in 2\mathbb{Z}_{>0}$.
For $\alpha,\beta\in\mathbb{R}_{>0}$ this difference equation has been proved
by Ismail and Rahman \cite{IR} using the explicit 
expression of $\mathcal{F}_\lambda(x)$ as very-well-poised ${}_8\phi_7$
series.
\end{proof}
Using the explicit expressions of $\mathcal{F}_\lambda(x)$, 
the conditions on $x,\lambda$
and the four parameters $\alpha,\rho,\beta,\sigma$ 
can be relaxed by meromorphic
continuation. The resulting function $\mathcal{F}_\lambda$ 
is a meromorphic, $\tau^{-1}$-periodic
eigenfunction of the Askey-Wilson second order 
difference operator $\mathcal{D}=
\mathcal{D}^{a,b,c,d}$ with
eigenvalue $E(\lambda)=E(\lambda;\beta)$.

\begin{rema} 
Some special cases of the matrix coefficients
$\phi_\lambda$ were explicitly expressed in terms of 
very-well-poised ${}_8\phi_7$ series in \cite{KS}
using the realization of the representation
$\pi_\lambda$ on the representation space $l^2(\mathbb{Z})$. In this approach
the basic hypergeometric series manipulations are much harder, since 
one needs a highly nontrivial evaluation of a non-symmetric Poisson type kernel
involving nonterminating ${}_2\phi_1$-series which is 
due to Rahman, see the appendix of
\cite{KS} and \cite{KR}. 
\end{rema}
\begin{rema}\label{Jacobirem}
The explicit ${}_8\phi_7$ expression \eqref{hs} of $\mathcal{F}_{\lambda}(x)$
and its analytic continuation was named the Askey-Wilson 
function in \cite{KS2}. 
Suslov \cite{S}, \cite{S2} established Fourier-Bessel type
orthogonality relations for the Askey-Wilson function.
Koelink and the author \cite{KS}, \cite{KS2} defined a generalized 
Fourier transform 
involving the Askey-Wilson function as the integral kernel, and
established its Plancherel and inversion formula. This transform,
called the Askey-Wilson function transform, arises as Fourier
transform on the noncompact quantum group $\hbox{SU}_q(1,1)$ (see \cite{KS}), 
and may thus be seen as a natural
analogue of the Jacobi function transform.
\end{rema}


\section{The expansion formula and the elliptic cosine kernel}

In this section we still assume that $\tau\in i\mathbb{R}_{>0}$, so 
$0<q=q_\tau=\exp(2\pi i\tau)<1$. To keep contact with the conventions
of the previous section, we keep working in base $q^2$.

First we recall the (normalized) Askey-Wilson polynomials \cite{AW}.
The Askey-Wilson polynomials are defined by
\[E_m(x)=E_m(x;a,b,c,d):={}_4\phi_3\left(\begin{matrix}
q^{-2m},q^{2m-2+a+b+c+d}, q^{a+x}, q^{a-x}\\
q^{a+b}, q^{a+c}, q^{a+d}\end{matrix}\,;\,q^2,q^2\right),\qquad 
m\in\mathbb{Z}_{\geq 0},
\]
with 
\[{}_4\phi_3\left(\begin{matrix} a_1,a_2,a_3,a_4\\
b_1,b_2,b_3\end{matrix}\,;\,q,z\right)=\sum_{k=0}^{\infty}
\frac{\bigl(a_1,a_2,a_3,a_4;q\bigr)_k}
{\bigl(q,b_1,b_2,b_3;q\bigr)_k}z^k,\qquad |z|<1.
\]
The Askey-Wilson polynomial $E_m(x)$ is a polynomial in $q^x+q^{-x}$
of degree $m$, normalized by $E_m(a)=1$. They
satisfy the orthogonality relations
\[\int_0^1E_m(x/\tau)E_n(x/\tau)\,
\frac{\bigl(q^{2x},q^{-2x};q^2\bigr)_{\infty}}
{\bigl(q^{a+x},q^{a-x},q^{b+x},q^{b-x},q^{c+x},q^{c-x},q^{d+x},q^{d-x};
q^2\bigr)_{\infty}}\,dx=0,\qquad m\not=n,
\]
provided that $\hbox{Re}(a),\hbox{Re}(b),\hbox{Re}(c),
\hbox{Re}(d)>0$. Observe that the Askey-Wilson polynomial 
$E_m(x)$ is regular at the special choice
\begin{equation}\label{choicepar}
(a,b,c,d)=\bigl(0,\frac{1}{2\tau},1,1-\frac{1}{2\tau}\bigr) 
\end{equation}
of Askey-Wilson parameters.
The above orthogonality relations also extend (by continuity) to
the Askey-Wilson parameters \eqref{choicepar}, leading to
\[\int_0^1E_m\bigl(x/\tau;0,\frac{1}{2\tau},1,1-\frac{1}{2\tau}\bigr)
E_n\bigl(x/\tau;0,\frac{1}{2\tau},1,1-\frac{1}{2\tau}\bigr)\,dx=0,
\qquad m\not=n,
\]
hence we conclude that
\begin{equation}\label{cosinepol}
E_m\bigl(x;0,\frac{1}{2\tau},1,1-\frac{1}{2\tau}\bigr)
=\cos(2\pi m\tau x),\qquad m\in\mathbb{Z}_{\geq 0},
\end{equation}
which is the usual cosine kernel from the Fourier theory on the unit circle.
In this section we derive an analogous result for the Askey-Wilson function.
                  
In the previous section we have introduced
the dual Askey-Wilson parameters associated to $a,b,c$ and $d$. 
They can be alternatively expressed as
\[(\widetilde{a},\widetilde{b},\widetilde{c},\widetilde{d})=
\bigl(-1+\frac{1}{2}(a+b+c+d), 1+\frac{1}{2}(a+b-c-d),
1+\frac{1}{2}(a-b+c-d),1+\frac{1}{2}(a-b-c+d)\bigr).
\]
Note that the special choice \eqref{choicepar} of Askey-Wilson parameters
is self dual,
$(a,b,c,d)=(\widetilde{a},\widetilde{b},\widetilde{c},\widetilde{d})$.
For our present purposes it is convenient to use yet another 
normalization of the Askey-Wilson function, namely
\begin{equation*}
\begin{split}
\mathfrak{E}^+(\mu,x)&=\mathfrak{E}^+(\mu,x;a,b,c,d)\\
&:=\frac{\bigl(q^{2+a-\widetilde{d}+\mu+x},q^{2+a-\widetilde{d}+\mu-x},
q^{2-a-d},q^{2+a-d};q^2\bigr)_{\infty}}
{\bigl(q^{\widetilde{a}+\widetilde{b}+\widetilde{c}+\mu},
q^{2-\widetilde{d}+\mu},q^{2-d+x},q^{2-d-x};q^2\bigr)_{\infty}}\\
&\quad\times{}_8W_7\bigl(
q^{-2+\widetilde{a}+\widetilde{b}+\widetilde{c}+\mu};
q^{a+x},q^{a-x},q^{\widetilde{a}+\mu},
q^{\widetilde{b}+\mu},q^{\widetilde{c}+\mu};
q^2,q^{2-\widetilde{d}-\mu}\bigr),\qquad |q^{2-\widetilde{d}-\mu}|<1.
\end{split}
\end{equation*}
For fixed $\lambda$, 
the eigenfunction $\mathcal{F}_\lambda(\cdot)$ of the 
Askey-Wilson second order difference operator $\mathcal{D}$ is 
a constant multiple of $\mathfrak{E}^+(2i\lambda,\cdot)$.
The present normalization of the Askey-Wilson function is convenient due
to the properties 
\begin{equation}\label{duality}
\mathcal{E}^+(\mu,x;a,b,c,d)=
\mathcal{E}^+(x,\mu;\widetilde{a},\widetilde{b},\widetilde{c},
\widetilde{d}\,)
\end{equation}
and $\mathfrak{E}^+(-\widetilde{a},-a)=1$. The property \eqref{duality}
is called {\it duality} and can be proved using 
a transformation formula for very-well-poised
${}_8\phi_7$ series, see \cite{KS2} for details. 
Furthermore,
\begin{equation}\label{polred}
\mathfrak{E}^+(\widetilde{a}+2m,x)=E_m(x),\qquad m\in \mathbb{Z}_{\geq 0},
\end{equation}
see e.g. \cite[(3.5)]{KS2}, thus the Askey-Wilson function 
$\mathfrak{E}^+(\mu,x)$ provides a natural meromorphic continuation of 
the Askey-Wilson polynomial in its degree.

The meromorphic continuation of  $\mathfrak{E}^+(\mu,x)$ in $\mu$ and $x$
can be established by the integral representation of the 
Askey-Wilson function 
(see the previous section), or by the expression of the Askey-Wilson function
as a sum of two balanced ${}_4\phi_3$'s (see e.g. \cite[(3.3)]{KS2}).
For our present purposes, it is most convenient to consider
the meromorphic continuation via the expansion formula 
of the Askey-Wilson function in Askey-Wilson polynomials, given by
\begin{equation}\label{ef}
\begin{split}
\mathfrak{E}^+(\mu,x)=&\frac{\bigl(q^{2-a-d},q^{2+a-d},q^{b+c},
q^{2+b-d},q^{2+c-d};q^2\bigr)_{\infty}}
{\bigl(q^{2+a+b+c-d},q^{2-d+x},q^{2-d-x},q^{2-\widetilde{d}+\mu},
q^{2-\widetilde{d}-\mu};q^2\bigr)_{\infty}}\\
\times&\sum_{m=0}^{\infty}E_m(x;a,b,c,2-d)E_m(\mu;\widetilde{a},\widetilde{b},
\widetilde{c},2-\widetilde{d}\,)\\
&\,\,\times 
\frac{\bigl(1-q^{2m+a+b+c-d}\bigr)\bigl(
q^{a+b+c-d},q^{a+b},q^{a+c};q^2\bigr)_m}
{\bigl(1-q^{a+b+c-d}\bigr)\bigl(
q^2,q^{2+b-d},q^{2+c-d};q^2\bigr)_m}(-1)^mq^{(1-a-d)m}q^{m^2},
\end{split}
\end{equation}
see \cite[Thm. 4.2]{St0}. 
The sum converges absolutely and uniformly on 
compacta of $(\mu,x)\in\mathbb{C}\times\mathbb{C}$ 
due to the Gaussian $q^{m^2}$. 
The expansion formula \eqref{ef} shows that
the Askey-Wilson function $\mathfrak{E}^+(\mu,x)$ is well defined and regular
at the special choice \eqref{choicepar} of Askey-Wilson parameters.
In fact, for this special choice of parameters, the Askey-Wilson
function can be expressed in terms of the (renormalized) Jacobi theta
function
\[\vartheta(x)=\bigl(-q^{1+x},-q^{1-x};q^2\bigr)_{\infty}
\]
as follows.
\begin{prop}\label{cosinethm}
We have the identity
\begin{equation}\label{cosine}
\mathfrak{E}^+\bigl(\mu,x;0,\frac{1}{2\tau},1,1-\frac{1}{2\tau}\bigr)=
\frac{\bigl(-q,-q;q^2\bigr)_{\infty}}{2}
\left(\frac{\vartheta(\mu+x)+\vartheta(\mu-x)}
{\vartheta(\mu)\vartheta(x)}\right).
\end{equation}
\end{prop}
\begin{proof}
To simplify notations, we write
\begin{equation}\label{E0}
\mathfrak{E}_0^+(\mu,x)=\mathfrak{E}^+\bigl(\mu,x;
0,\frac{1}{2\tau},1,1-\frac{1}{2\tau}\bigr)
\end{equation}
for the duration of the proof.

We substitute the special choice \eqref{choicepar} of Askey-Wilson
parameters in the expansion formula for $\mathfrak{E}^+(\mu,x)$.
By simple $q$-series manipulations and by \eqref{cosinepol}, 
we obtain the explicit formula
\[\mathfrak{E}_0^+(\mu,x)=
\frac{\bigl(-q,-q,-q,q,-q^2;q^2\bigr)_{\infty}}
{\bigl(q^2;q^2\bigr)_{\infty}\vartheta(\mu)\vartheta(x)}
\left(1+2\sum_{m=1}^{\infty}\cos(2\pi m\tau\mu)\cos(2\pi m\tau x)\,q^{m^2}
\right).
\]
Using the well known Jacobi triple product identity 
\[\sum_{m=-\infty}^{\infty}q^{m^2+mx}=
1+2\sum_{m=1}^{\infty}\cos(2\pi m\tau x)\,q^{m^2}=
\bigl(q^2;q^2\bigr)_{\infty}\,
\vartheta(x),
\]
and the elementary identity
\[\cos(2\pi m\tau\mu)\cos(2\pi m\tau x)=\frac{1}{2}\bigl(
\cos(2\pi m\tau (\mu+x))+\cos(2\pi m\tau(\mu-x))\bigr),
\]
we deduce that
\[\mathfrak{E}_0^+(\mu,x)=\frac{\bigl(-q,-q,-q,q,-q^2;q^2\bigr)_{\infty}}
{2}\left(\frac{\vartheta(\mu+x)+\vartheta(\mu-x)}
{\vartheta(\mu)\vartheta(x)}\right).
\]
Simplifying the multiplicative constant yields the desired result.
\end{proof}

\begin{rema}\label{class}
The Jacobi theta function $\vartheta(x)$ is the natural
$\tau^{-1}$-periodic analogue of the Gaussian $q^{-x^2}$, since
\[\sum_{k=-\infty}^{\infty}q^{-(x+k\tau^{-1})^2}=\sqrt{-2i\tau}\,
\bigl(q^2;q^2\bigr)_{\infty}
\vartheta(x)
\]
by the Jacobi triple product identity and the Jacobi inversion formula.
In fact, in \cite{St0} and \cite{St} it is shown that
the function $\bigl(q^{2-d+x},q^{2-d-x};q^2\bigr)_{\infty}$,
which reduces to $\vartheta(x)$ for the Askey-Wilson parameters 
\eqref{choicepar}, plays the role of the Gaussian in the Askey-Wilson theory.
If one replaces the theta functions by Gaussians 
in the right hand side of \eqref{cosine}, then
we obtain up to a multiplicative constant
\[\frac{1}{2}\left(\frac{q^{-(\mu+x)^2}+q^{-(\mu-x)^2}}{q^{-\mu^2-x^2}}\right)=
\frac{1}{2}\bigl(q^{-2\mu x}+q^{2\mu x}\bigr)=\cos(4\pi\tau\mu x),\]
which is essentially the classical cosine kernel. Thus the
right hand side of \eqref{cosine} is an elliptic analogue of
the cosine kernel.
\end{rema}
\begin{rema}
Using the quasi-periodicity 
\begin{equation}\label{qp}
\vartheta(x+2)=q^{-1-x}\vartheta(x)
\end{equation}
of the Jacobi theta function, we obtain as
a consequence of \eqref{cosine},
\[\mathfrak{E}^+\bigl(2m,x;0,\frac{1}{2\tau},1,1-\frac{1}{2\tau}\bigr)=
\cos(2\pi m\tau x),\qquad m\in\mathbb{Z}_{\geq 0},
\]
which is in accordance with \eqref{cosinepol} and \eqref{polred}.
\end{rema}
\begin{rema}
The orthogonality relations for the Askey-Wilson polynomials
with Askey-Wilson parameters \eqref{choicepar}
are equivalent to the $L^2$-theory of 
the classical Fourier tansform on the unit circle. 
On the other hand, the $L^2$-theory of the Askey-Wilson function
transform, see \cite{KS2}, does {\it not} reduce to the
$L^2$-theory of the classical Fourier theory on the real line
for the Askey-Wilson parameters \eqref{choicepar}. 
Instead one obtains a Fourier type transform with integral kernel given by
the elliptic cosine function \eqref{cosine}. For the corresponding 
$L^2$ theory, the transform is defined on a weighted 
$L^2$-space consisting of functions that are supported on a 
finite closed interval and an infinite, unbounded sequence of 
discrete mass points. 
This transform, as well as the general Askey-Wilson
function transform, still has many properties in common with 
the classical Fourier transform on the real line, 
see e.g. \cite{KS2}, \cite{St0} and \cite{St}.
\end{rema}

Cherednik's \cite{C} Hecke algebra approach to $q$-special functions
leads to a direct proof that the right hand side 
of the expansion formula \eqref{ef} is 
an eigenfunction of the Askey-Wilson second order difference 
operator $\mathcal{D}$,
see \cite{St}. 
The expansion formula \eqref{ef} may thus be seen as the 
explicit link between Cherednik's approach and Ismail's and Rahman's 
\cite{IR} construction of eigenfunctions of $\mathcal{D}$ in terms of
very-well-poised ${}_8\phi_7$ series. 
We end this section by sketching a proof of
Proposition \ref{cosinethm} using Cherednik's Hecke algebra approach. 

The affine Hecke algebra techniques for Askey-Wilson polynomials 
are developed in full detail in \cite{NS}, and for
Askey-Wilson functions in \cite{St}. 
We first recall one of the main results from \cite{St}, specialized to the
present rank one situation.

Define two difference-reflection operators by
\begin{equation}\label{T}
\begin{split}
\bigl(T_0^{c,d}f\bigr)(x)&=-q^{-2+c+d}f(x)+\frac{(1-q^{c-x})(1-q^{d-x})}
{(1-q^{2-2x})}\,\bigl(f(2-x)-f(x)\bigr),\\
\bigl(T_1^{a,b}f\bigr)(x)&=-q^{a+b}f(x)+
\frac{(1-q^{a+x})(1-q^{b+x})}{(1-q^{2x})}\,\bigl(f(-x)-f(x)\bigr).
\end{split}
\end{equation}
The connection with affine 
Hecke algebras follows from the fact that $T_0=T_0^{c,d}$
and $T_1=T_1^{a,b}$ satisfy
Hecke type quadratic relations. These relations imply
that the operators $T_0$ and $T_1$ are invertible.
Consider the (invertible) operator
\begin{equation}\label{Yn}
Y=Y^{a,b,c,d}:=T_1^{a,b}\circ T_0^{c,d}.
\end{equation}
\begin{rema}\label{linkd}
The operator $Y+Y^{-1}$, acting on even functions, is essentially the
Askey-Wilson second order difference operator $\mathcal{D}$, see e.g. 
\cite[Prop. 5.8]{NS}. 
\end{rema}
Theorem 5.17 in \cite{St} states that for generic Askey-Wilson parameters
$(a,b,c,d)$,  
there exists a unique meromorphic function $\mathfrak{E}(\cdot,\cdot)=
\mathfrak{E}(\cdot,\cdot;a,b,c,d)$
on $\mathbb{C}\times\mathbb{C}$ 
satisfying the following six conditions:
\begin{enumerate}
\item[{\bf 1.}] $\mathfrak{E}(\mu,x)$ is $\tau^{-1}$-periodic in $\mu$ and $x$,
\item[{\bf 2.}] $(\mu,x)\mapsto \bigl(q^{2-\widetilde{d}+\mu},
q^{2-\widetilde{d}-\mu},q^{2-d+x},
q^{2-d-x};q^2\bigr)_{\infty}\,\mathfrak{E}(\mu,x)$
is analytic,
\item[{\bf 3.}] For fixed generic $\mu\in\mathbb{C}$, 
$\mathfrak{E}(\mu,\cdot)$ is an eigenfunction
of $Y^{a,b,c,d}$ with eigenvalue $q^{\widetilde{a}-\mu}$,
\item[{\bf 4.}] For fixed generic $x\in\mathbb{C}$, $\mathfrak{E}(\cdot,x)$ is
an eigenfunction of 
$Y^{\widetilde{a},\widetilde{b},\widetilde{c},\widetilde{d}}$
with eigenvalue $q^{a-x}$,
\item[{\bf 5.}] $\bigl(T_1^{a,b}\mathfrak{E}(\mu,\cdot)\bigr)(x)=
-q^{a+b}\,
\bigl(T_1^{\widetilde{a},\widetilde{b}}\mathfrak{E}(\cdot,x)\bigr)(\mu)$,
\item[{\bf 6.}] $\mathfrak{E}(-\widetilde{a},-a)=1$.
\end{enumerate}
The existence of a kernel $\mathfrak{E}$ satisfying the above six conditions
is proved by explicitly constructing $\mathfrak{E}$ as series
expansion in nonsymmetric analogues of the
Askey-Wilson polynomials, see \cite[(6.6)]{St}. 
This expansion formula for $\mathfrak{E}$
is very similar to the expansion 
formula \eqref{ef} of the Askey-Wilson function
$\mathfrak{E}^+$ in Askey-Wilson polynomials. 
In fact, a comparison of the
formulas leads to the explicit link
\begin{equation}\label{connection}
\mathfrak{E}^+(\mu,x)=\bigl(C^+_{a,b}\mathfrak{E}(\mu,\cdot)\bigr)(x),\qquad
C^+_{a,b}:=\frac{1}{1-q^{a+b}}\bigl(1+T_1^{a,b}\bigr),
\end{equation}
see \cite[Thm. 6.20]{St}. These results allow us to study
the Askey-Wilson function $\mathfrak{E}^+$ using the
characterizing conditions {\bf 1--6} for the underlying kernel $\mathfrak{E}$,
instead of focussing on the explicit expression for $\mathfrak{E}^+$.

The kernel $\mathfrak{E}(\mu,x)$ is regular at the 
Askey-Wilson parameters \eqref{choicepar}. The resulting kernel
\[\mathfrak{E}_0(\mu,x)=\mathfrak{E}\bigl(\mu,x;0,\frac{1}{2\tau},
1,1-\frac{1}{2\tau}\bigr),
\]
is the unique meromorphic kernel satisfying the six conditions 
{\bf 1--6} for
the special Askey-Wilson parameters \eqref{choicepar}. Observe that
the operators $T_0,T_1$ and $Y$ for the special Askey-Wilson parameters
\eqref{choicepar} reduce to
\[(T_0f)(x)=f(2-x),\qquad (T_1f)(x)=f(-x),\qquad (Yf)(x)=f(2+x),
\]
hence $\mathfrak{E}_0$ is the unique 
meromorphic kernel satisfying the six conditions
\begin{enumerate}
\item[{\bf 1${}^\prime$.}] $\mathfrak{E}_0(\mu,x)$ is 
$\tau^{-1}$-periodic in $\mu$ and $x$,
\item[{\bf 2${}^\prime$.}] 
$(\mu,x)\mapsto \vartheta(\mu)\vartheta(x)\mathfrak{E}_0(\mu,x)$
is analytic,
\item[{\bf 3${}^\prime$.}] 
$\mathfrak{E}_0(\mu,x+2)=q^{-\mu}\mathfrak{E}_0(\mu,x)$,
\item[{\bf 4${}^\prime$.}] 
$\mathfrak{E}_0(\mu+2,x)=q^{-x}\mathfrak{E}_0(\mu,x)$,
\item[{\bf 5${}^\prime$.}] $\mathfrak{E}_0(\mu,-x)=\mathfrak{E}_0(-\mu,x)$,
\item[{\bf 6${}^\prime$.}] $\mathfrak{E}_0(0,0)=1$.
\end{enumerate}
We conclude that 
\[\mathfrak{E}_0(\mu,x)=
\bigl(-q,-q;q^2\bigr)_{\infty}\,
\frac{\vartheta(\mu+x)}{\vartheta(\mu)\vartheta(x)},
\]
since the right hand side satisfies {\bf 1${}^\prime$--6${}^\prime$}
due to the quasi-periodicity \eqref{qp} of $\vartheta(x)$. 
Thus $\mathfrak{E}_0(\mu,x)$ is an elliptic analogue 
of the exponential kernel $\exp(-4\pi i\tau\mu x)$, cf.  Remark \ref{class}.
Using the notation \eqref{E0}, we conclude that
\[\mathfrak{E}_0^+(\mu,x)=
\bigl(C^+_{0,\frac{1}{2\tau}}\mathfrak{E}_0(\mu,\cdot)
\bigr)(x)=\frac{1}{2}\bigl(\mathfrak{E}_0(\mu,x)+\mathfrak{E}_0(\mu,-x)\bigr),
\]
which is the desired formula \eqref{cosine}.


\section{The Askey-Wilson function for $|q|=1$.}

In this section we take $-\frac{1}{2}<\tau<0$, so that
$q=q_\tau=\exp(2\pi i\tau)$ has modulus one and $q\not=\pm 1$. 
The assignment
\begin{equation}
(K^{\pm 1})^\star=K^{\pm 1},\qquad (X^{\pm })^\star=-X^{\pm}
\end{equation}
uniquely extends to a unital, anti-linear, anti-algebra involution
on $\mathcal{U}_q$. This particular choice of $\star$-structure 
corresponds to the real form $\mathfrak{s}\mathfrak{l}(2,\mathbb{R})$ of
$\mathfrak{s}\mathfrak{l}(2,\mathbb{C})$. The $\star$-unitary sesquilinear
form for the representations $\pi_\lambda$ ($\lambda\in\mathbb{R}$)
of $\mathcal{U}_q$ (see Section 3) is
\[\langle f,g\rangle^\prime=\int_{-i\infty}^{i\infty}f(z)\overline{g(z)}\,dz,
\]
where $f$ and $g$ are meromorphic functions 
which are regular on a large enough strip around the 
imaginary axes and decay sufficiently fast 
at $\pm i\infty$. Koornwinder's twisted primitive
element $iY_\rho\in \mathcal{U}_q$ is $\star$-selfadjoint for 
$\rho\in \mathbb{R}$. Thus in principle we are 
all set to extend the construction 
of eigenfunctions of the Askey-Wilson second order difference operator
$\mathcal{D}$ to the $|q|=1$ case by simply replacing the role
of the $q$-gamma function $\Gamma_{2\tau}$ by $G_{2\tau}$. 
We need to be careful though due to the following differences with the 
$0<q<1$ case:
\begin{enumerate}
\item[{\bf a.}] The analogue of the
explicit eigenfunction of $iY_\rho$ (cf. \eqref{f}), given now as
quotient of hyperbolic gamma functions $G_{2\tau}$, has more singularities.
\item[{\bf b.}] No $\tau^{-1}$-periodicity conditions have to be imposed. 
Consequently, the parameters $\alpha$ and $\beta$ do not need to
be discretized for the $|q|=1$ case.
\item[{\bf c.}] We have to take the decay rates at $\pm i\infty$
of integrands into account.
\end{enumerate}
One needs to be careful with the
decay rate (see {\bf c}) in reproving the crucial Theorem 
\ref{main1} because acting by $\pi_\lambda(X^\pm)$
worsens the asymptotics at $\pm i\infty$ 
(the factor $q^{z}=\exp(2\pi i\tau z)$
is $\mathcal{O}(\exp(-2\pi\tau\hbox{Im}(z))$ as 
$\hbox{Im}(z)\rightarrow \infty$). The decay rate can be improved by
considering different eigenfunctions of $iY_\rho$, but then the
location of the singularities turns out to cause problems.

To get around these problems, we generalize the techniques
of Section 4 to a {\it nonunitary} set-up. More concretely, we replace
$\langle \cdot,\cdot\rangle^\prime$ by a bilinear form, given as
a contour integral over a certain deformation of $i\mathbb{R}$.
With such a bilinear form,  
the singularity problems and the asymptotic problems 
can be resolved simultaneously.

The proper
replacement of the involution $\star$ is the unique unital, {\it linear},
anti-algebra involution $\circ$ on $\mathcal{U}_q$ satisfying
\[\bigl(K^{\pm 1}\bigr)^\circ=K^{\mp 1},\qquad 
\bigl(X^{\pm}\bigr)^\circ=-X^{\pm}.
\]
For $f\in \mathcal{M}$ we write $S(f)\subset \mathbb{C}$ for the
singular set of $f$. 
\begin{defi}
$f\in \mathcal{M}$ is said to have \textup{(}exponential\textup{)}
growth rate $\epsilon\in\mathbb{R}$
at $\pm i\infty$ when the following two conditions are satisfied:

{\bf 1.} For some compact subset $K_f\subset \mathbb{R}$,
\[S(f)\subset \{z\in \mathbb{C} \, | \, \hbox{Im}(z)\in K_f\}.
\]

{\bf 2.} The function $f$ satisfies
\[|f(x+iy)|=\mathcal{O}\bigl(\exp(\epsilon|y|)\bigr),\qquad y\rightarrow
\pm \infty,
\]
uniformly for $x$ in compacts of $\mathbb{R}$.
\end{defi}
We call a contour $\mathcal{C}$ a {\it deformation of} $i\mathbb{R}$
when $\mathcal{C}$ intersects the line
\[l_c=\{z\in \mathbb{C} \, | \, \hbox{Im}(z)=c\} 
\]
in exactly one point $z_c(\mathcal{C})$ for all $c\in \mathbb{R}$,
and $z_c(\mathcal{C})=ic$ for $|c|>>0$. 
For $k\in\mathbb{Z}_{\geq 0}$ we define the strip of radius $k$ around
$\mathcal{C}$ by
\[\bigcup_{c\in\mathbb{R}}\{z\in l_c \,\, | \,\, |z-z_c(\mathcal{C})|\leq k \}.
\] 
For $k=0$, the strip of radius zero around $\mathcal{C}$ is the contour
$\mathcal{C}$ itself.
\begin{lem}\label{altsa}
Suppose that $f,g\in \mathcal{M}$ have growth rates $\epsilon_f$
and $\epsilon_g$ at $\pm i\infty$, respectively. Suppose furthermore that
$\epsilon_f+\epsilon_g<2\pi\tau$ and 
that $f$ and $g$ are analytic 
on the strip of radius one around a given \textup{(}oriented\textup{)}
 deformation $\mathcal{C}$ of $i\mathbb{R}$. Then
\[\bigl(\pi_\lambda(X)f,g\bigr)_{\mathcal{C}}=
\bigl(f,\pi_{-\lambda}(X^\circ)g\bigr)_{\mathcal{C}},
\qquad \forall\, X\in\mathcal{U}_q^1,
\]
where
\[\bigl(f,g\bigr)_{\mathcal{C}}=\int_\mathcal{C}f(z)g(z)dz.
\]
\end{lem} 
\begin{proof}
The proof follows by an elementary application of Cauchy's Theorem.
\end{proof}
\begin{rema}\label{better}
The condition on the growth rates in Lemma \ref{altsa} may be weakened
to $\epsilon_f+\epsilon_g<0$ 
when $X$ is taken from the subspace 
$\hbox{span}_{\mathbb{C}}\{1,K^{-1},K\}$ of $\mathcal{U}_q^1$.
\end{rema}

Koornwinder's twisted primitive element $Y_\rho$ (see \eqref{Y}) is not
$\circ$-invariant,
\[Y_\rho^\circ=-q^{\frac{1}{2}}K^{-1}X^++q^{-\frac{1}{2}}K^{-1}X^-+
\left(\frac{q^{-\rho}+q^{\rho}}{q^{-1}-q}\right)(K^{-2}-1).
\]
On the other hand, a direct computation shows that 
$\pi_{-\lambda}(Y_\rho^\circ)$ is still a first order difference operator.
This leads to the following analogue of Proposition \ref{diffone}. 
\begin{prop}\label{difftwo}
Let $\alpha,\rho\in\mathbb{C}$.
A meromorphic function $f\in\mathcal{M}$ is an eigenfunction of
$\pi_{-\lambda}(Y_{\rho}^\circ)$ with eigenvalue $\mu_\alpha(\rho)$
if and only if 
\[f(z-2)=\frac{\sin(\pi(i\lambda+\alpha+\rho-z)\tau)
\sin(\pi(-i\lambda+\alpha+\rho+z)\tau)}
{\sin(\pi(i\lambda+\rho-1+z)\tau)
\sin(\pi(-i\lambda+\rho+1-z)\tau)}\,f(z).
\]
\end{prop}
We define two meromorphic functions by
\begin{equation}\label{newform}
\begin{split}
g_\lambda(z;\beta,\sigma)&=
\frac{G_{2\tau}\bigl(\frac{1}{2\tau}-1+\beta+\sigma-i\lambda+z\bigr)
G_{2\tau}\bigl(\frac{1}{2\tau}+\sigma-i\lambda-z\bigr)}
{G_{2\tau}\bigl(-\frac{1}{2\tau}+1+\beta+\sigma+i\lambda-z\bigr)
G_{2\tau}\bigl(-\frac{1}{2\tau}+\sigma+i\lambda+z\bigr)},\\
h_\lambda(z;\alpha,\rho)&=
\frac{G_{2\tau}\bigl(-\frac{1}{2\tau}-1+\alpha+\rho+i\lambda-z\bigr)
G_{2\tau}\bigl(-\frac{1}{2\tau}+\rho+i\lambda+z\bigr)}
{G_{2\tau}\bigl(-\frac{1}{2\tau}+1+\alpha+\rho-i\lambda+z\bigr)
G_{2\tau}\bigl(-\frac{1}{2\tau}+\rho-i\lambda-z\bigr)}.
\end{split}
\end{equation}
The difference equation for $G_{2\tau}$, Proposition \ref{diffone}
and Proposition \ref{difftwo} imply
\begin{equation}\label{efalt}
\pi_\lambda(Y_\sigma)g_\lambda(\cdot;\beta,\sigma)=
\mu_\beta(\sigma)g_\lambda(\cdot;\beta,\sigma),\qquad
\pi_{-\lambda}(Y_\rho^\circ)h_\lambda(\cdot;\alpha,\rho)=
\mu_\alpha(\rho)h_\lambda(\cdot;\alpha,\rho).
\end{equation}
We want to construct now eigenfunctions of the gauged Askey-Wilson second
order difference operator $\mathcal{L}$ which are of the form 
\begin{equation}\label{psi}
\begin{split}
\psi_\lambda(x;\alpha,\rho,\beta,\sigma)&=
\bigl(\pi_\lambda(\widehat{x}K)g_\lambda(\cdot;\beta,\sigma),
h_\lambda(\cdot;\alpha,\rho)\bigr)_{\mathcal{C}}\\
&=\int_{\mathcal{C}}g_\lambda(1+x+z;\beta,\sigma)h_\lambda(z;\alpha,\rho)dz
\end{split}
\end{equation}
for a suitable deformation $\mathcal{C}$ of $i\mathbb{R}$. To make
sense of this integral, we need to take the singularities and the 
asymptotic behaviour at $\pm i\infty$
of the integrand into account. The singularities can be located
using the precise information on the zeros and poles of the hyperbolic
gamma function $G_\tau$, see Proposition \ref{Ruij}{\bf (ii)}.
It follows that the singularities of 
$z\mapsto g_\lambda(1+x+z;\beta,\sigma)$ are contained
in the union of the four half lines
\begin{equation*}
\begin{split}
&-1-\beta-\sigma+i\lambda-x+\mathbb{R}_{\leq 0},\qquad
-1+\beta+\sigma+i\lambda-x+\mathbb{R}_{\leq 0},\\
&-\sigma-i\lambda-x+\mathbb{R}_{\geq 0},\qquad\qquad\qquad
\sigma-i\lambda-x+\mathbb{R}_{\geq 0},
\end{split}
\end{equation*}
and the singularities of $z\mapsto h_\lambda(z;\alpha,\rho)$ 
are contained in the
union of the four half lines
\begin{equation*}
\begin{split}
&-1+\frac{1}{\tau}-\rho-i\lambda+\mathbb{R}_{\leq 0},\qquad\quad\,\,\,\,
-1+\rho-i\lambda+\mathbb{R}_{\leq 0},\\
&-\alpha-\rho+i\lambda+\mathbb{R}_{\geq 0},\qquad\qquad\quad\,\,\,
-\frac{1}{\tau}+\alpha+\rho+i\lambda+\mathbb{R}_{\geq 0}.
\end{split}
\end{equation*}
We call the above eight half lines the {\it singular half lines}
with respect to the given, fixed parameters 
$\tau, x,\lambda,\alpha,\rho,\beta,\sigma$.
Each singular half line is contained in some horizontal
line $l_c$ for some $c\in \mathbb{R}$. For generic 
parameters $\alpha,\rho,\beta,\sigma$, the eight singular
half lines lie on different horizontal lines.
Under these generic assumptions, there exists a deformation 
$\mathcal{C}_x$ of $i\mathbb{R}$ which
seperates the four singular half lines with real part 
tending to $-\infty$ from the 
four singular half lines with real part tending to $\infty$.
We take such contour $\mathcal{C}=\mathcal{C}_x$ 
in the definition \eqref{psi} of $\psi_\lambda$.
The resulting function
$\psi_\lambda(x)$ is well defined and independent of 
the particular choice of the deformation $\mathcal{C}_x$ of 
$i\mathbb{R}$, since 
$g_\lambda(\cdot;\beta,\sigma)$
(respectively $h_\lambda(\cdot;\alpha,\rho)$)
has growth rate $\pi\bigl((1-2\hbox{Im}(\lambda))\tau-1\bigr)$
(respectively $\pi\bigl(1+2\hbox{Im}(\lambda)\bigr)\tau$)
at $\pm i\infty$. This follows
from the asymptotic behaviour of the hyperbolic gamma function
$G_\tau$, see Proposition \ref{Ruij}{\bf (iv)}. The resulting function
$\psi_\lambda(x)$ is analytic in $x$.
\begin{thm}
For generic parameters $\alpha,\rho,\beta,\sigma$, the 
function $\psi_\lambda(\cdot)=\psi_\lambda(\cdot;\alpha,\rho,\beta,\sigma)$
defined by
\begin{equation}\label{psinew}
\begin{split}
\psi_\lambda(x)&=\int_{\mathcal{C}_x}g_\lambda(1+x+z;\beta,\sigma)
h_\lambda(z;\alpha,\rho)\,dz\\
&=\int_{\mathcal{C}_x}
\frac{G_{2\tau}\bigl(\frac{1}{2\tau}+\beta+\sigma-i\lambda+x+z\bigr)
G_{2\tau}\bigl(\frac{1}{2\tau}-1+\sigma-i\lambda-x-z\bigr)}
{G_{2\tau}\bigl(-\frac{1}{2\tau}+\beta+\sigma+i\lambda-x-z\bigr)
G_{2\tau}\bigl(-\frac{1}{2\tau}+1+\sigma+i\lambda+x+z\bigr)}\\
&\qquad\qquad\quad
\times \frac{G_{2\tau}\bigl(-\frac{1}{2\tau}-1+\alpha+\rho+i\lambda-z\bigr)
G_{2\tau}\bigl(-\frac{1}{2\tau}+\rho+i\lambda+z\bigr)}
{G_{2\tau}\bigl(-\frac{1}{2\tau}+1+\alpha+\rho-i\lambda+z\bigr)
G_{2\tau}\bigl(-\frac{1}{2\tau}+\rho-i\lambda-z\bigr)}\,dz
\end{split}
\end{equation}
is an eigenfunction of the gauged Askey-Wilson difference operator
$\mathcal{L}^{\alpha,\rho,\beta,\sigma}$ with eigenvalue $E(\lambda;\beta)$.
\end{thm}
\begin{proof}
We adjust the proof of Theorem \ref{main1} to the present set-up.
We simplify notations by writing $g(z)=g_\lambda(z;\beta,\sigma)$
and $h(z)=h_\lambda(z;\alpha,\rho)$.
Choose the deformation $\mathcal{C}_x$ of
$i\mathbb{R}$ such that $g$ and $h$ are analytic on the strip of radius 
$\geq 4$ around $\mathcal{C}_x$. 

Observe that 
$\pi_\lambda(\widehat{x}\Omega K)g$ has the same growth rate at $\pm i\infty$
as $g$. By \eqref{Omegaaction} and the definition of 
$\psi_\lambda(x)$, we conclude that the integral
$\bigl(\pi_\lambda(\widehat{x}\Omega K)g,h\bigr)_{\mathcal{C}_x}$ converges
absolutely and equals
\[
\left(\frac{q^{i\lambda}-q^{-i\lambda}}{q-q^{-1}}\right)^2\,\psi_\lambda(x).
\]
On the other hand, 
the radial part computation of $\Omega$ with respect to Koornwinder's twisted
primitive element yields
\begin{equation}\label{radeq}
\widehat{x}\Omega K=\widehat{x}\Omega(x)K+
\bigl((Y_\rho-\mu_\alpha(\rho))K^{-1}\bigr)X\widehat{x}K 
+\widehat{x}Z\bigl(Y_\sigma-\mu_\beta(\sigma)\bigr)
\end{equation}
for certain elements $X,Z\in\mathcal{U}_q^1$, cf. Proposition \ref{radial}.
Substituting this algebraic identity in 
$\bigl(\pi_\lambda(\widehat{x}\Omega K)g,h\bigr)_{\mathcal{C}_x}$
and using $\pi_\lambda(Y_\sigma)g=\mu_\beta(\sigma)g$, we have
\begin{equation}\label{two}
\bigl(\pi_\lambda(\widehat{x}\Omega K)g,h\bigr)_{\mathcal{C}_x}
=\bigl(\pi_\lambda(\widehat{x}\Omega(x)K)g,h\bigr)_{\mathcal{C}_x}
+\bigl(\pi_\lambda\bigl((Y_\rho-\mu_\alpha(\rho))K^{-1}X\widehat{x}K\bigr)g,
h\bigr)_{\mathcal{C}_x},
\end{equation}
provided that both integrals on the right hand side of \eqref{two}
converge absolutely.

For the second term on the right hand side of \eqref{two}, observe that
the sum of the growth rates of $g$ and $h$ at $\pm i\infty$
equals $(2\tau-1)\pi$. Furthermore,
\[(2\tau-1)\pi<4\pi\tau
\]
since $-\frac{1}{2}<\tau<0$, and
\[(Y_\rho-\mu_\alpha(\rho))K^{-1}\in\mathcal{U}_q^1.
\]
Hence the integral 
\[\bigl(\pi_\lambda\bigl((Y_\rho-\mu_\alpha(\rho))K^{-1}X\widehat{x}K\bigr)g,
h\bigr)_{\mathcal{C}_x}
\]
converges absolutely, and Lemma \ref{altsa} shows that this integral
equals
\[\bigl(\pi_\lambda(X\widehat{x}K)g,\pi_{-\lambda}(K)
\pi_{-\lambda}(Y_\rho^\circ-\mu_\alpha(\rho))h\bigr)_{\mathcal{C}_x}=0,
\]
where the last equality follows from the eigenvalue equation
$\pi_{-\lambda}(Y_\rho^\circ)h=\mu_\alpha(\rho)h$.

For the first term on the right hand side of 
\eqref{two}, observe that  
$\pi_\lambda(\widehat{x}\Omega(x)K)g$ has the same growth rate at
$\pm i\infty$ as $g$. 
Hence the integral $\bigl(\pi_\lambda\bigl(\widehat{x}\Omega(x)K\bigr)g,
h\bigr)_{\mathcal{C}_x}$ converges absolutely.
By the definition of the gauged Askey-Wilson second order difference
operator $\mathcal{L}$, this integral equals
\[\frac{q^{\beta-1}}{(q-q^{-1})^2}\bigl(\mathcal{L}\psi_\lambda\bigr)(x)
+\frac{q^{\beta-1}(1-q^{1-\beta})^2}
{(q-q^{-1})^2}\,\psi_\lambda(x).
\]

Combining the results, we conclude that
\[\bigl(\pi_\lambda(\widehat{x}\Omega K)g,h\bigr)_{\mathcal{C}_x}=
\bigl(\pi_\lambda(\widehat{x}\Omega(x)K)g,h\bigr)_{\mathcal{C}_x},
\]
and for both sides we have obtained an explicit expression in terms of
$\psi_\lambda$. The resulting identity for $\psi_\lambda$
yields the desired difference equation
$\mathcal{L}\psi_\lambda=E(\lambda)\psi_\lambda$.
\end{proof}
\begin{rema}
The eigenfunction $\psi_\lambda$ of $\mathcal{L}$ as given
by the integral \eqref{psinew} looks very similar to the eigenfunction
$\phi_\lambda$ of $\mathcal{L}$ (for the case $\tau\in i\mathbb{R}_{>0}$)
as given by the integral \eqref{phi2}, since the integrand of $\phi_\lambda$
is essentially the integrand of $\psi_\lambda$ with
the hyperbolic gamma functions $G_{2\tau}$ replaced
by $q$-gamma functions $\Gamma_{2\tau}$. 
Note though that the integration cycles are different. 
\end{rema}
We can reformulate the difference equation for $\psi_\lambda$
in terms of the Askey-Wilson
second order difference operator $\mathcal{D}=\mathcal{D}^{a,b,c,d}$
(see \eqref{D}) using an appropriate gauge factor, cf. Corollary \ref{D1}.
Using the parameter correspondence \eqref{parlinks}, we can for instance
choose the gauge factor $\delta(x)=\delta(x;\alpha,\rho,\beta,\sigma)$
by
\begin{equation}\label{delta}
\delta(x)=\frac{G_{2\tau}\bigl(-1+\frac{1}{2\tau}+a-x\bigr)
G_{2\tau}\bigl(-1+\frac{1}{2\tau}+b-x\bigr)}
{G_{2\tau}\bigl(1+\frac{1}{2\tau}-c-x\bigr)
G_{2\tau}\bigl(1+\frac{1}{2\tau}-d-x\bigr)}.
\end{equation}
\begin{cor}
For generic parameters $\alpha,\rho,\beta,\sigma$, the function 
$H_\lambda=H_\lambda(\cdot;a,b,c,d)\in\mathcal{M}$ defined by
$H_\lambda(x)=\delta(x)^{-1}\psi_\lambda(x)$ 
satisfies the Askey-Wilson second order difference equation
\[(\mathcal{D}H_\lambda)(x)=E(\lambda)H_\lambda(x),
\] 
where the Askey-Wilson parameters $(a,b,c,d)$
are given by \eqref{parlinks}. 
\end{cor}

\section{Appendix}
In this appendix we sketch a proof of 
Proposition \ref{radial}, following closely
the arguments of Koornwinder \cite{K} for the special case $\alpha=\beta=0$.

Fix $x\in\mathbb{C}$ and
write $X\equiv X^\prime$ for $X,X^\prime
\in\widetilde{\mathcal{U}}_q$ if
\[X-X^\prime\in (Y_{\rho}-\mu_\alpha(\rho))\,\widehat{x}\,\mathcal{U}_q^1+
\widehat{x}\,\mathcal{U}_q^1(Y_{\sigma}-\mu_\beta(\sigma)).
\]
In order to simplify notations, I use the notations
$\mu=\mu_\alpha(\rho)$ and $\nu=\mu_\beta(\sigma)$ for the duration of the
proof. To reduce $\widehat{x}\,\Omega K=\widehat{x}\,K\Omega$ in the desired form, 
we need to concentrate on the part
$\widehat{x}\,KX^+X^-$. Using
\begin{equation*}
\begin{split}
\widehat{x}\,KX^+X^-=&q^{\frac{1}{2}+x}\bigl(q^{\frac{1}{2}}X^+K-
q^{-\frac{1}{2}}X^-K\bigr)\widehat{x}X^-\\
&+
q^{-\frac{1}{2}+x}X^-\widehat{x}\bigl(q^{-\frac{1}{2}}X^-K-
q^{\frac{1}{2}}X^+K\bigr)+q^{2x}\widehat{x}KX^-X^+
\end{split}
\end{equation*}
and the commutation relation between $X^-$ and $X^+$ in $\mathcal{U}_q$,
we obtain
\begin{equation*}
\begin{split}
(1-q^{2x})\,\widehat{x}\,KX^+X^-\equiv
&q^{\frac{1}{2}+x}
\left(\frac{q^\rho+q^{-\rho}}{q-q^{-1}}(K^2-1)+\mu\right)\widehat{x}\,X^-\\
&-q^{-\frac{1}{2}+x}X^-\widehat{x}\,
\left(\frac{q^\sigma+q^{-\sigma}}{q-q^{-1}}(K^2-1)+\nu\right)
+q^{2x}\left(\frac{K^{-1}-K^3}{q-q^{-1}}\right)\widehat{x},
\end{split}
\end{equation*}
hence we need to focus now only on the reduction of $\widehat{x}X^-$
and $\widehat{x}K^2X^-$. Using
\begin{equation*}
\begin{split}
\widehat{x}\,X^-=
&q^{\frac{1}{2}-x}
\bigl(q^{-\frac{1}{2}}X^-K-q^{\frac{1}{2}}X^+K\bigr)K^{-1}\widehat{x}\\
&+q^{\frac{3}{2}-2x}K^{-1}\widehat{x}\,
\bigl(q^{\frac{1}{2}}X^+K-q^{-\frac{1}{2}}X^-K\bigr)
+q^{1-2x}K^{-1}\widehat{x}\,X^-K,
\end{split}
\end{equation*}
we obtain
\begin{equation*}
\begin{split}
(1-q^{2-2x})\,\widehat{x}\,X^-\equiv
&-q^{\frac{1}{2}-x}\left(
\frac{q^\rho+q^{-\rho}}{q-q^{-1}}(K^2-1)+\mu\right)K^{-1}\widehat{x}\\
&\qquad+q^{\frac{3}{2}-2x}K^{-1}\widehat{x}\,
\left(\frac{q^\sigma+q^{-\sigma}}{q-q^{-1}}(K^2-1)+\nu\right),
\end{split}
\end{equation*}
and a similar computation yields
\begin{equation*}
\begin{split}
(1-q^{-2-2x})\,\widehat{x}\,K^2X^-\equiv
&-q^{-\frac{3}{2}-x}\left(
\frac{q^\rho+q^{-\rho}}{q-q^{-1}}(K^2-1)+\mu\right)K\widehat{x}\\
&\qquad+q^{-\frac{5}{2}-2x}K\widehat{x}\,
\left(\frac{q^\sigma+q^{-\sigma}}{q-q^{-1}}(K^2-1)+\nu\right).
\end{split}
\end{equation*}
Collecting all these results we arrive for generic 
$x\in\mathbb{C}$ at a formula of the form
\[\widehat{x}\,\Omega K\equiv \widehat{x}\,
(\widetilde{B}(x)K^3+\widetilde{C}(x)K+\widetilde{D}(x)K^{-1})
\]
for explicit rational expressions
$\widetilde{B}(x)$, $\widetilde{C}(x)$ and $\widetilde{D}(x)$.
Using
\begin{equation*}
\begin{split}
q^\rho+q^{-\rho}-(q-q^{-1})\mu&=q^{\rho+\alpha}+q^{-\rho-\alpha},\\
q^\sigma+q^{-\sigma}-(q-q^{-1})\nu&=q^{\sigma+\beta}+q^{-\sigma-\beta},
\end{split}
\end{equation*}
the rational functions $\widetilde{B}(x)$, $\widetilde{C}(x)$ and 
$\widetilde{D}(x)$ are explicitly given by
\begin{equation*}
\begin{split}
\widetilde{B}(x)=&\frac{\bigl(q^{\frac{1}{2}+x}(q^\rho+q^{-\rho})-
q^{\frac{3}{2}+2x}(q^\sigma+q^{-\sigma})\bigr)
\bigl(q^{-\frac{5}{2}-2x}(q^\sigma+q^{-\sigma})-
q^{-\frac{3}{2}-x}(q^\rho+q^{-\rho})\bigr)}
{(q-q^{-1})^2(1-q^{2x})(1-q^{-2-2x})}\\
&\qquad-\frac{q^{2x}}{(q-q^{-1})(1-q^{2x})}
+\frac{q^{-1}}{(q-q^{-1})^2},
\end{split}
\end{equation*}
\begin{equation*}
\begin{split}
&\widetilde{C}(x)=-\frac{2}{(q-q^{-1})^2}\\
&+\frac{\bigl(q^{\frac{1}{2}+x}(q^\rho+q^{-\rho})-
q^{\frac{3}{2}+2x}(q^\sigma+q^{-\sigma})\bigr)
\bigl(q^{-\frac{3}{2}-x}(q^{\rho+\alpha}+q^{-\rho-\alpha})-
q^{-\frac{5}{2}-2x}(q^{\sigma+\beta}+q^{-\sigma-\beta})\bigr)}
{(q-q^{-1})^2(1-q^{2x})(1-q^{-2-2x})}\\
&+\frac{\bigl(q^{-\frac{1}{2}+2x}(q^{\sigma+\beta}+q^{-\sigma-\beta})
-q^{\frac{1}{2}+x}(q^{\rho+\alpha}+q^{-\rho-\alpha})\bigr)
\bigl(q^{\frac{3}{2}-2x}(q^\sigma+q^{-\sigma})
-q^{\frac{1}{2}-x}(q^\rho+q^{-\rho})\bigr)}
{(q-q^{-1})^2(1-q^{2x})(1-q^{2-2x})},
\end{split}
\end{equation*}
and
\begin{equation*}
\begin{split}
\widetilde{D}(x)&=
\frac{q^{2x}}{(q-q^{-1})(1-q^{2x})}+\frac{q}{(q-q^{-1})^2}\\
&+\bigl(q^{-\frac{1}{2}+2x}(q^{\sigma+\beta}+q^{-\sigma-\beta})
-q^{\frac{1}{2}+x}(q^{\rho+\alpha}+q^{-\rho-\alpha})\bigr)\\
&\qquad\times
\frac{\bigl(q^{\frac{1}{2}-x}(q^{\rho+\alpha}+
q^{-\rho-\alpha})-q^{\frac{3}{2}-2x}
(q^{\sigma+\beta}+q^{-\sigma-\beta})\bigr)}
{(q-q^{-1})^2(1-q^{2x})(1-q^{2-2x})}.
\end{split}
\end{equation*}
By a  tedious computation, it can now
be proven that
\begin{equation*}
\begin{split}
\widetilde{B}(x)&=\frac{q^{\beta-1}}{(q-q^{-1})^2}B(x),\\
\widetilde{C}(x)&=\frac{q^{\beta-1}}{(q-q^{-1})^2}
\bigl(-A(x)-A(x^{-1})+(1-q^{1-\beta})^2\bigr),\\
\widetilde{D}(x)&=\frac{q^{\beta-1}}{(q-q^{-1})^2}D(x),
\end{split}
\end{equation*}
with the parameters $a,b,c,d$ as in \eqref{parlinks}.

\end{document}